\documentclass[11pt]{article}
\pdfoutput=1

\usepackage[backend=biber,
            bibstyle=phys,
            sorting=none,
            url=false,
            doi=false,
            eprint=true,
            sortcites=true,
            citestyle=numeric-comp]{biblatex}

\DeclareFieldFormat*{title}{\mkbibitalic{#1}}
\DeclareFieldFormat{labelnumberwidth}{\mkbibbrackets{#1}~~}

\usepackage[paperwidth=8.5in, paperheight=11in, portrait, top=1in, bottom=1.in,
 left=1.15in, right=1.15in]{geometry}
\usepackage[utf8]{inputenc}
\usepackage{amssymb,amsthm,amsmath,amsfonts,mathtools,bm,bbold}
\usepackage{graphicx,units,mathrsfs,latexsym,hyperref,xcolor,subfigure,epsfig}
\usepackage{tikz} \usetikzlibrary{calc,decorations.markings,positioning,shapes,snakes}

\def\sl{{\mathfrak{sl}}}
\def\cQ{{\mathcal{Q}}}

\def\saw{\mathbf{saw}}
\def\zh{{Zh_q}}
\def\szh{{sZh_q}}
\def\Hq{\widehat{H}_q}
\def\SHq{S\!{H}_q}
\newcommand{\cR}{\mathcal{R}}
\newcommand{\czR}{\check{\mathcal{R}}}
\newcommand{\mC}{\mathfrak{C}}

\DeclareMathOperator{\id}{id}
\DeclareMathOperator{\Sk}{Sk}

\def\noi{\noindent}

\numberwithin{equation}{section}
\theoremstyle{plain}
\newtheorem{thm}{Theorem}[section]

\newtheorem{prop}[thm]{Proposition}
\newtheorem{coro}[thm]{Corollary}

\newtheorem{defi}[thm]{Definition}

\newtheorem{rmk}[thm]{Remark}
\newtheorem{rem}[thm]{Remark}

\newcommand{\KP}[1]{\begin{tikzpicture}
 [baseline=-\dimexpr\fontdimen22\textfont2\relax] #1 \end{tikzpicture}}
\newcommand{\KPA}{\KP{\filldraw[color=blue, fill=none, thick] circle (0.3);} }
\newcommand{\KPB}{\KP{\draw[color=blue,thick] (-0.3,0.3) -- (0.3,-0.3);
 \draw[color=blue,thick] (-0.3,-0.3) -- (-0.05,-0.05);
 \draw[color=blue,thick] (0.05,0.05) -- (0.3,0.3);} }
\newcommand{\KPC}{\KP{\draw[color=blue,thick](-0.3,0.3) .. controls (0,-0.05)..(0.3,0.3);
 \draw[color=blue,thick](-0.3,-0.3) .. controls (0,0.05) .. (0.3,-0.3);} }
\newcommand{\KPD}{\KP{\draw[color=blue,thick](-0.3,-0.3)..controls (0.05,0)..(-0.3,0.3);
 \draw[color=blue,thick] (0.3,-0.3) .. controls (-0.05,0) .. (0.3,0.3);} }
\newcommand{\punctiii}{\draw (0,0) node{$\times$}; \draw (1,0) node{$\times$};
 \draw (2,0) node{$\times$};} 
\newcommand{\punctiv}{\draw (0,0) node{$\times$}; \draw (1,0) node{$\times$};
 \draw (2,0) node{$\times$}; \draw (3,0) node{$\times$};} 
\newcommand{\larc}[1]{\draw (#1,-0.25cm) arc (270:90:0.25cm);
} 
\newcommand{\rarc}[1]{\draw (#1,-0.25cm) arc (-90:90:0.25cm);
} 
\newcommand{\uarc}[1]{\draw (#1 cm-0.25cm,0) arc (180:0:0.25cm);
} 
\newcommand{\darc}[1]{\draw (#1 cm-0.25cm,0) arc (180:360:0.25cm);
} 
\newcommand{\updown}[2]{\draw (#1,-0.25cm) -- (#2,-0.25cm);
\draw (#1,0.25cm) -- (#2,0.25cm);} 
\newcommand{\arrul}[1]{\draw (#1+1,0.35cm) arc (60:120:1.cm);
\draw (#1,0.35cm) -- (#1+0.125,0.525cm); \draw (#1,0.35cm) -- (#1+0.2,0.35cm);
} 
\newcommand{\arrur}[1]{\draw (#1+1,0.35cm) arc (60:120:1.cm);
\draw (#1+1,0.35cm) -- (#1+0.875,0.525cm); \draw (#1+1,0.35cm) -- (#1+0.8,0.35cm);
} 
\newcommand{\arrdr}[1]{\draw (#1,-0.35cm) arc (-120:-60:1cm);
\draw (#1+1,-0.35cm) -- (#1+0.875,-0.525cm); \draw (#1+1,-0.35cm) -- (#1+0.8,-0.35cm);
} 
\newcommand{\arrdl}[1]{\draw (#1,-0.35cm) arc (-120:-60:1cm);
\draw (#1,-0.35cm) -- (#1+0.125,-0.525cm); \draw (#1,-0.35cm) -- (#1+0.2,-0.35cm);
} 
\newcommand{\bendu}[2]{\draw[color=black] (#1 cm,0) to [bend left=90] (#2 cm,0);}
\newcommand{\bendd}[2]{\draw[color=black] (#1 cm,0) to [bend right=90] (#2 cm,0);}
\newcommand{\punct}[1]{\draw (#1 cm,0) node{$\times$};}

\begin{document}
\title{\bf The Askey--Wilson algebra and its avatars}
\author{
Nicolas Crampé\textsuperscript{$1$}
\footnote{E-mail: crampe1977@gmail.com}~,
Luc Frappat\textsuperscript{$2$}
\footnote{E-mail: luc.frappat@lapth.cnrs.fr}~,
Julien Gaboriaud\textsuperscript{$3$}
\footnote{E-mail: julien.gaboriaud@umontreal.ca}~,
Loïc Poulain d'Andecy\textsuperscript{$4$}
\footnote{E-mail: loic.poulain-dandecy@univ-reims.fr}~,\\
Eric Ragoucy\textsuperscript{$2$}
\footnote{E-mail: eric.ragoucy@lapth.cnrs.fr}~,
Luc Vinet\textsuperscript{$3$}
\footnote{E-mail: vinet@CRM.UMontreal.CA}~,
\\[.5em]
\textsuperscript{$1$}
\small Institut Denis-Poisson CNRS/UMR 7013 - Université de Tours - Université
d'Orléans,\\
\small~Parc de Grandmont, 37200 Tours, France.\\[.9em]
\textsuperscript{$2$}
\small Laboratoire d'Annecy-le-Vieux de Physique Th\'eorique LAPTh,\\
\small~Université Grenoble Alpes, Université Savoie Mont Blanc, CNRS, F-74000 Annecy,
 France.\\[.9em]
\textsuperscript{$3$}
\small Centre de Recherches Math\'ematiques, Universit\'e de Montr\'eal, \\
\small~P.O. Box 6128, Centre-ville Station, Montr\'eal (Qu\'ebec), H3C 3J7,
Canada.\\[.9em]
\textsuperscript{$4$}
\small Laboratoire de mathématiques de Reims UMR 9008,
Université de Reims Champagne-Ardenne,\\
\small~Moulin de la Housse BP 1039, 51100 Reims, France.
}
\date{\today} \maketitle

\noindent{\bf Abstract:}
The original Askey--Wilson algebra introduced by Zhedanov encodes the bispectrality
properties of the eponym polynomials.  The name \textit{Askey--Wilson algebra} is
currently used to refer to a variety of related structures that appear in a large number
of contexts. We review these versions, sort them out and establish the relations between
them.  We focus on two specific avatars. The first is a quotient of the original Zhedanov
algebra; it is shown to be invariant under the Weyl group of type $D_4$
and to have a reflection algebra presentation.  The second is a universal analogue of the
first one; it is isomorphic to the Kauffman bracket skein algebra (KBSA) of the
four-punctured sphere and to a subalgebra of the universal double affine Hecke algebra
$(C_1^{\vee},C_1)$. This second algebra emerges from the Racah problem of $U_q(\sl_2)$ and
is related via an injective homomorphism to the centralizer of $U_q(\sl_2)$ in its
threefold tensor product. How the Artin braid group acts on the incarnations of this
second avatar through conjugation by $R$-matrices (in the Racah problem) or half Dehn
twists (in the diagrammatic KBSA picture) is also highlighted.  Attempts at defining
higher rank Askey--Wilson algebras are briefly discussed and summarized in a diagrammatic
fashion.
\\[.5em]

\clearpage

\tableofcontents\par

\clearpage

\section{Introduction}\label{sec:intro}


In order to provide an algebraic underpinning for the Askey--Wilson polynomials
\cite{AskeyWilson1985}, Zhedanov introduced what he called the Askey--Wilson algebra
\cite{Zhedanov1991}. We shall refer to it rather as the Zhedanov algebra. The
Askey--Wilson polynomials sit at the top of the Askey classification scheme of the
hypergeometric orthogonal polynomials \cite{KoekoekLeskyetal2010} and are, consequently,
of fundamental interest; their algebraic interpretation by Zhedanov hence bears
commensurate importance.  These $q$-polynomials are bispectral: in addition to verifying a
three-term recurrence prescribed by Favard's theorem for any family of orthogonal
polynomials \cite{Chihara2011}, they are also eigenfunctions of a $q$-difference operator.
The Zhedanov algebra was constructed by taking these two bispectral operators as
generators and identifying the relations they obey. As sometimes happens with natural
constructs, related structures have emerged in a variety of contexts and have typically
all been called Askey--Wilson algebras. This propensity keeps rising and it is hence
timely to review the topic.  This paper will provide a taxonomy and a description of the
algebras that loosely go under the name of \textit{Askey--Wilson algebras} and  will
characterize in some depth two avatars of particular relevance. It will also set the stage
for the exploration of generalizations.

The focus of this survey will be on algebraic aspects. Before we discuss the contents in
more details, let us briefly go over some of the manifestations of these Askey--Wilson
algebras and the advances they have generated. Grosso modo, they have had direct
applications in physical models and have also been at the heart of mathematical
developments establishing useful interconnections between fields.  One occurrence is in
the recoupling of three irreducible representations of $U_q(\mathfrak{sl}_2)$ which is
called its Racah problem.  It is known that the $6j$-symbols of this algebra are expressed
in terms of $q$-Racah polynomials which are a finite truncation of the Askey--Wilson ones.
As a rule, whenever the Askey--Wilson polynomials (or their truncated version) appear, the
associated algebra will be present. In the case of the Racah problem, it is found that the
intermediate Casimir elements verify Askey--Wilson relations \cite{GranovskiiZhedanov1993,
Huang2016}.  These polynomials and algebras appear in the study of the ASEP model with
open boundaries \cite{UchiyamaSasamotoetal2004}, as martingale polynomials and quadratic
harnesses in probabilistic models \cite{BrycMatysiaketal2007} and are connected to (a
degeneration of) the Sklyanin algebra \cite{GorskyZabrodin1993, WiegmannZabrodin1995,
GaboriaudTsujimotoetal2020}. Quite generally, the Askey--Wilson algebras are present in
the context of integrable models, through the Yang--Baxter and reflection equations
\cite{Baseilhac2005, Baseilhac2005a, BaseilhacKoizumi2005, Aneva2008,
AnevaChaichianetal2008, VinetZhedanov2008}, and can be viewed as truncations of the
$q$-Onsager algebra \cite{Baseilhac2005}.  Elements of representation theory have been
investigated in \cite{Zhedanov1991, NoumiStokman2004, Huang2015, Huang2016, Huang2020} and
another of its manifestations is as a coideal subalgebra of $U_q(\sl_2)$
\cite{GranovskiiZhedanov1993a, Terwilliger2011a, CrampeShaabankabakibo2021}. The
Askey--Wilson algebras have also been cast in the framework of Howe duality using the pair
$(U_q(\sl_2), \mathfrak{o}_{q^{1/2}}(2n))$ \cite{FrappatGaboriaudetal2020,
GaboriaudVinetetal2020, FrappatGaboriaudetal2019, NoumiUmedaetal1996}; they are special
cases of the recently introduced Painlev\'e algebras \cite{Mazzocco2016} and belong to the
Calabi--Yau class \cite{EtingofGinzburg2010}.  There is a significant connection to the
field of algebraic combinatorics, as Askey--Wilson algebras are central in the
classification of $P$- and $Q$- polynomial association schemes and the study of Leonard
pairs and triples \cite{BannaiIto1984,Terwilliger2001, TerwilligerVidunas2004,
Terwilliger2011, Huang2012, Terwilliger2018}.  The Askey--Wilson algebras have also been
shown to offer a promising platform to extend the quantum Schur--Weyl duality to arbitrary
representations and have been seen in that respect to admit the Temperley--Lieb  and
Birman--Murakami--Wenzl algebras \cite{CrampeVinetetal2020a} as quotients.  Askey--Wilson
algebras have moreover found their way in the general framework of knot theory through
their identification with the Kauffman bracket skein algebras of the four-punctured sphere
$\Sk_{\,iq^{1/2}}(\Sigma_{0,4})$ and other elementary surfaces \cite{BullockPrzytycki1999,
Cooke2020, CookeLacabanne2020}.  This is also closely connected to double affine Hecke
algebras (DAHA) as the Askey--Wilson algebra is related to the spherical subalgebra of the
DAHA of type $(C_1^\vee,C_1)$ \cite{Oblomkov2004, Koornwinder2007a, Koornwinder2008,
ItoTerwilliger2010, Terwilliger2013, Mazzocco2016, KoornwinderMazzocco2018,
TsujimotoVinetetal2019a, Huang2020}.

This overview of the relevance of Askey--Wilson algebras in different domains motivates
the present topical report.  Let us make at this point a few additional remarks on the
introduction of the algebra $\Sk_{\,iq^{1/2}}(\Sigma_{0,4})$ in the Askey--Wilson picture
to stress that this paper also features novel results relating the Askey--Wilson algebra,
the Kauffman bracket skein algebra and the braid group.

Kauffman bracket skein algebras (KBSA) have been defined independently by Turaev
\cite{Turaev1991} and Bullock and Przytycki \cite{BullockPrzytycki1999} in the study of
knot invariants and can be seen to encompass the celebrated Jones polynomial
\cite{Jones1985, Kauffman1987}.  Computations in the KBSA are done through diagrammatic
manipulations given by a set of rules (the skein relations). It is appreciated that this
$\Sk_{\,iq^{1/2}}(\Sigma_{0,4})$ algebra is closely related to the centralizer of
$U_q(\sl_2)$ in its threefold tensor product. This ties in with the Temperley--Lieb
algebra which admits a diagrammatic presentation \cite{Kauffman1987, Kauffman1990,
KauffmanLins1994} for generic $q$, is precisely the centralizer of $U_q(\sl_2)$ in the
threefold tensor product of the fundamental representations of $U_q(\sl_2)$
\cite{Jimbo1985} and, as already indicated, was found to be a quotient of the
Askey--Wilson algebra \cite{CrampeVinetetal2020a}.

A natural question that has arisen asks about higher rank extensions of Askey--Wilson
algebras. In view of the ubiquity of the 3-generated Askey--Wilson algebras it is to be
expected that such generalizations will prove quite fruitful. This question is non-trivial
however since many avenues that are likely to yield different outcomes can be followed.
Among those possibilities, one is to consider the algebra realized by the intermediate
Casimir elements in multifold tensor products of $U_q(\sl_2)$ \cite{PostWalter2017,
DeBieDeClercqetal2020, DeBieDeClercq2020, DeClercq2019}, and another is to increase the
rank of the algebra $U_q(\sl_2)$ to, say, $U_q(\sl_3)$ when studying the Racah problem.
Augmenting the number of punctures of the sphere in the KBSA approach could also be
envisaged.  Making much sense is the idea to start from the multivariate Askey--Wilson
polynomials \cite{GasperRahman2005}, to work out the algebra formed by its bispectral
operators \cite{Iliev2011, BaseilhacMartin2018, DeBieDeClercq2020} and to take things from
there. This is after all how the story began.  Steps have been taken in these directions
but final conclusions have not been reached. Some authors have considered higher order
truncations of the reflection algebra \cite{BaseilhacKoizumi2005a} understood as a
quotient of the $q$-Onsager algebra (see also \cite{BaseilhacCrampe2019} for the classical
limit of this result).  The upshot is that there is currently no clear consensus on what
the higher rank Askey--Wilson algebra is\footnote{Remarkably, for the $q\to1$ and $q\to-1$
limits of the Askey--Wilson algebra, higher rank extensions have been more successfully
defined respectively in \cite{DeBieGenestetal2017} for the Racah algebra and in
\cite{DeBieGenestetal2016} for the Bannai--Ito algebra.}. This is not too surprising since
there are still a few loose ends in the rank one cases.

As a prelude to a solid understanding of the higher rank Askey--Wilson algebra, it is
appropriate to clarify the picture for the ordinary Askey--Wilson algebras. Indeed, as
these algebras have appeared in multiple instances in the literature, names, conventions
and notations are quite diverse. We are here proposing a standardization and offering a
number of new results.  The paper will unfold as follows.  The various Askey--Wilson
avatars  will be introduced in Section \ref{sec:aws}. They will be given names and defined
in a comparative way. Emphasis will be put on two particular versions.
The first is a
quotient of the Zhedanov algebra which we will call the Special Zhedanov algebra. In
Section \ref{sec:reflection}, we will show that the Zhedanov algebra is obtained as the
reflection algebra defined from particular $R$- and reflection matrices.  In this
formalism, the Special Zhedanov algebra corresponds to fixing the Sklyanin determinant to
a certain value; the name \textit{Special} is chosen in analogy with the nomenclature of
Lie groups.  A Weyl group $W(D_4)$ symmetry of the Special Zhedanov algebra will then be
presented in Section \ref{sec:szh}, thus generalizing an analoguous result for the Racah
algebra.
The second avatar that will be closely looked at will be called the Special
Askey--Wilson algebra. It can be seen as the equivalent of the Special Zhedanov algebra
where the parameters are promoted to central elements in the algebra.  That this algebra
is isomorphic to the Kauffman bracket skein algebra of the four-punctured sphere
$\Sk_{\,iq^{1/2}}(\Sigma_{0,4})$ is the object of Section \ref{sec:KBS}. In Section
\ref{sec:cent}, the Special Askey--Wilson algebra will further be related to the algebra
$\mathcal{A}_3$ associated to the Racah problem of $U_q(\sl_2)$ and to the centralizer
$\mC_3$ of $U_q(\sl_2)$ in its threefold tensor product. An injective homomorphism of
algebras between the latter two structures will be stated and its proof will be found in
Appendix \ref{app:inj}.  The relation between the Special Askey--Wilson algebra and the
universal double affine Hecke algebra (DAHA) of type $(C_1^{\vee},C_1)$ will be discussed
in Section \ref{sec:daha}. How the Artin braid group $B_3$ acts on both the
$\mathcal{A}_3$ and $\Sk_{\,iq^{1/2}}(\Sigma_{0,4})$ algebras, respectively through
conjugation by braided $R$-matrices and through half Dehn twists will be highlighted in
Section \ref{sec:braid}.
The question of the possible higher-rank generalizations of the
Askey--Wilson algebra will be addressed in Section \ref{sec:higher}. A crossing index will
be introduced and used to summarize efficiently the main results of \cite{PostWalter2017}
and \cite{DeClercq2019} and new relations for the higher rank analogues will be provided.
Elements of interest for further study of the higher rank generalizations of the Special
Askey--Wilson algebra will be offered in addition. Concluding remarks will end the paper.

\section{Askey--Wilson algebras}\label{sec:aws}

\subsection{A jungle of Askey--Wilson algebras}
\label{sec:aw3}

As mentioned in the above, the name \textit{Askey--Wilson algebra} has appeared and been
connected to diverse objects in a multitude of contexts. Therefore, the notations and
appellations in the literature are sometimes confusing. For the sake of clarity, we start
by presenting these different algebraic structures and give to them unambiguous names to
distinguish them.

\paragraph{The Askey--Wilson algebra $\mathbf{aw}(3)$} is the unital associative
algebra depending on the parameter $q$ with generators $C_{12}$, $C_{23}$, $C_{13}$ and
central elements $C_1$, $C_2$, $C_3$, $C_{123}$ obeying the $\mathbb{Z}_3$-symmetric
relations
\begin{subequations}\label{eq:aw3}
\begin{align}
 \label{eq:aw3A}
 C_{12}+\frac{[C_{23},C_{13}]_q}{q^{2}-q^{-2}}&=\frac{C_1C_2+C_3C_{123}}{q+q^{-1}},\\
 \label{eq:aw3B}
 C_{23}+\frac{[C_{13},C_{12}]_q}{q^{2}-q^{-2}}&=\frac{C_2C_3+C_1C_{123}}{q+q^{-1}},\\
 \label{eq:aw3C}
 C_{13}+\frac{[C_{12},C_{23}]_q}{q^{2}-q^{-2}}&=\frac{C_3C_1+C_2C_{123}}{q+q^{-1}},
\end{align}
where the $q$-commutator is defined by $[A,B]_q=qAB-q^{-1}BA$.
Throughout the paper, we suppose that $q\in\mathbb{C}$ is not a root of unity.
The Casimir element of this algebra is
\begin{align}\label{eq:aw3cas}
\begin{aligned}
 \varOmega:=qC_{12}C_{23}C_{13}+q^{2}{C_{12}}^{2}+q^{-2}{C_{23}}^{2}+q^{2}{C_{13}}^{2}
 &-qC_{12}(C_1C_2+C_3C_{123})\\
 -q^{-1}C_{23}(C_2C_3+C_1C_{123})&-qC_{13}(C_3C_1+C_2C_{123}).
\end{aligned}
\end{align}
\end{subequations}
{Let us emphasize that this algebra $\mathbf{aw}(3)$ is not the algebra called
Askey--Wilson algebra by A. Zhedanov, and denoted $AW(3)$ in \cite{Zhedanov1991}. In the
present paper, we call the latter the \textit{Zhedanov algebra} (see below).\\[1em]
\noindent From the $\mathbf{aw}(3)$ algebra, we define multiple quotients or subalgebras
which appear in different contexts; these justify the importance of this algebra.}

\paragraph{The Special Askey--Wilson algebra $\saw(3)$} is the quotient of
$\mathbf{aw}(3)$ by the supplementary relation
\begin{align}\label{eq:awt3cas1}
 \varOmega&=(q+q^{-1})^{2}-{C_{123}}^{2}-{C_{1}}^{2}-{C_{2}}^{2}-{C_{3}}^{2}
 -C_{123}C_{1}C_{2}C_{3}.
\end{align}
A justification of the adjective \textit{special} is given in Section
\ref{sec:reflection}.  This algebra is isomorphic to the Kauffman bracket skein module of
the four-punctured sphere (see Section \ref{sec:KBS}) and is directly associated to the
centralizer of the diagonal action of $U_q(\sl_2)$ in its threefold tensor product (see
Section \ref{sec:cent}).

\paragraph{The universal Askey--Wilson algebra $\Delta_q$} defined in
\cite{Terwilliger2011} is the subalgebra of $\mathbf{aw}(3)$ generated by $C_{12}$,
$C_{23}$, $C_{13}$ as well as the central elements $\alpha=C_1C_2+C_3C_{123}$,
$\beta=C_2C_3+C_1C_{123}$ and $\gamma=C_3C_1+C_2C_{123}$.  The Casimir element of
$\Delta_q$ becomes
\begin{align}\label{eq:AW3PTcas}
 \Omega=qC_{12}C_{23}C_{13}+q^{2}C_{12}^{2}+q^{-2}C_{23}^{2}+q^{2}C_{13}^{2}
 -qC_{12}\alpha-q^{-1}C_{23}\beta-qC_{13}\gamma.
\end{align}
{An injective homomorphism of $\Delta_q$ into $U_q(\sl_2)\otimes U_q(\sl_2)\otimes
U_q(\sl_2)$ has been studied in \cite{Huang2016} and
its finite irreducible representations have been classified in \cite{Huang2015}.
The universal Askey--Wilson algebra also intersects the theory of free Lie algebras, see
e.g. \cite{Cantuba2015} and \cite{CantubaSilvestrov2020}.}

\paragraph{The evaluated Askey--Wilson algebra $Z_q(m_1,m_2,m_3)$} is the
quotient of $\mathbf{aw}(3)$ by the supplementary relations
\begin{align}
 C_i=q^{m_i}+q^{-m_i}, \qquad i=1,2,3.
\end{align}
It plays a central role in the study of the centralizer of the diagonal embedding of
$U_q(\sl_2)$ in the threefold tensor product of representations of $U_q(\sl_2)$
\cite{CrampeVinetetal2020a}.

\paragraph{The Zhedanov algebra $\zh(m_1,m_2,m_3,m_4)$} is the quotient of
$\mathbf{aw}(3)$ by
\begin{align}\label{eq:mm}
 C_i=q^{m_i}+q^{-m_i}, \qquad C_{123}=q^{m_4}+q^{-m_4}, \qquad i=1,2,3,
\end{align}
and was first introduced by Zhedanov as the algebra encoding the bispectrality of the
Askey--Wilson polynomials \cite{Zhedanov1991}. To be precise, in \cite{Zhedanov1991}, an
alternative equivalent presentation recalled in \eqref{eq:K1}--\eqref{eq:K3},
\eqref{eq:xxi2}--\eqref{eq:xxi4p} has been given.
The above $\mathbb{Z}_3$-symmetric presentation of $\zh(m_1,m_2,m_3,m_4)$ is introduced in
\cite{WiegmannZabrodin1995}. {This algebra appears to be also the proper algebraic
setting to characterize the Leonard pairs \cite{TerwilligerVidunas2004}.}

\paragraph{The Special Zhedanov algebra $\szh(m_1,m_2,m_3,m_4)$} is obtained as the
quotient of $\saw(3)$ by relations \eqref{eq:mm} (see \eqref{eq:K1}--\eqref{eq:xxi6} for
an alternative presentation).  It appears naturally as the commutation relations of the
intermediate Casimir elements acting on the multiplicity space of the decomposition of the
threefold tensor product of representations of $U_q(\sl_2)$ (see Section \ref{sec:Zmm}).

\subsection{Miscellaneous properties}\label{sec:pro}

\subsubsection*{PBW basis}

The Askey--Wilson algebra $\mathbf{aw}(3)$ has a Poincar\'e--Birkhoff--Witt (PBW) basis
given explicitly by the following elements
\begin{align}
 C_{12}^{i} C_{23}^{j} C_{13}^{k} C_1^{m} C_2^{n} C_3^{p} C_{123}^{q},\qquad
 \qquad i,j,k,m,n,p,q\in\mathbb{N}.
\end{align}
The proof is a slight generalization of the proof of the PBW basis for the universal Askey--Wilson algebra $\Delta_q$
given in \cite{Terwilliger2011}.
We can also obtain a PBW basis for the Special Askey--Wilson algebra $\saw(3)$ from the
one of $\mathbf{aw}(3)$ by restricting the range of the exponent $j$ to $\{0,1\}$ instead
of $\mathbb{N}$.

\subsubsection*{Calabi--Yau algebra}

The Zhedanov algebra $\zh(m_1,m_2,m_3,m_4)$ can be derived from a Calabi--Yau potential in
the following sense \cite{Ginzburg2007}.  Let $F=\mathbb{C}[x_1, x_2, x_3]$ be a free
associative algebra and view $F$ as a graded algebra such that $deg(x_1)=d_1$,
$deg(x_2)=d_2$ and $deg(x_3)=d_3$ (with $0<d_1\leq d_2 \leq d_3$).  We define
$F_{cycl}=F/[F,F]$ and the map $\tfrac{\partial}{\partial x_j}: F_{cycl} \rightarrow F$ on
cyclic words as follows
\begin{equation}
\frac{\partial [x_{i_1}x_{i_2}\dots x_{i_r}] }{\partial x_j}
 = \sum_{\{s|i_s=j\}} x_{i_s +1}x_{i_s +2}\dots x_{i_r}x_{i_1}x_{i_2}\dots x_{i_s -1}
\end{equation}
and we extend it to $F_{cycl}$ by linearity.  Let $\Phi(x_1,x_2, x_3) \in F_{cycl}$ be a
potential which can be decomposed as follows
\begin{equation}
 \Phi(x_1,x_2, x_3)=\Phi^{(d)}(x_1,x_2, x_3)+\Phi^{<d}(x_1,x_2, x_3),
\end{equation}
where $\Phi^{(d)}(x_1,x_2, x_3)$ is homogeneous of degree $d=d_1+d_2+d_3$ and
$\Phi^{<d}(x_1,x_2, x_3)$ is composed of terms of degree strictly inferior to $d$.  Then
the algebra whose defining relations are given by
\begin{equation}\label{eq:dphi}
\frac{\partial\Phi}{\partial x_j}=0, \qquad j=1,2,3,
\end{equation}
is a Calabi--Yau algebra \cite{EtingofGinzburg2010}.

Now, let $x_1=K_{12}$, $x_2=K_{23}$, $x_3=K_{13}$ and $deg(x_1)=deg(x_2)=2$, $deg(x_3)=3$.
Consider the potential
\begin{align}\label{eq:Phi}
\begin{aligned}
\Phi^{(7)}(x_1,x_2, x_3)&=q[x_1x_2x_3]-q^{-1}[x_1x_3x_2],\\
\Phi^{<7}(x_1,x_2, x_3)&=(q+q^{-1})( [x_1x_2^2]+[x_1^2 x_2])-\xi_4 [x_1] -\xi'_4 [x_2]
 -\tfrac{1}{2} [x_3^2] -\xi_2 [x_1x_2] .
\end{aligned}
\end{align}
It is easy to see that the defining relations of $\zh(m_1,m_2,m_3,m_4)$ presented in
\eqref{eq:K1}-\eqref{eq:K3} are equivalent to imposing \eqref{eq:dphi} for the potential
\eqref{eq:Phi}.  In other words, $\zh(m_1,m_2,m_3,m_4)$ derives from the Calabi--Yau
potential $\Phi$ \eqref{eq:Phi}.

\section{The Zhedanov algebra as a truncated reflection algebra}\label{sec:reflection}

In this section, we recall \cite{Baseilhac2005} that the defining relations of the algebra
$\zh(m_1,m_2,m_3,m_4)$ can be equivalently encoded in a reflection equation
\cite{Sklyanin1988}.  This realization of an algebra is usually called the FRT
presentation, in honor of the authors of \cite{FaddeevReshetikhinetal1990}.  This
presentation allows one to connect the Zhedanov algebra to the reflection algebra which is
intensively studied in the context of quantum integrable systems.  In addition, we show
that the algebra $\szh(m_1,m_2,m_3,m_4)$ can be also obtained naturally by setting the
Sklyanin determinant to a certain value; this justifies the appellation \textit{special}
for the quotiented algebra since it is obtained by fixing the value of a determinant, as
in the definition of the Special Linear group $SL_n$.

The cornerstone of the FRT presentation is the $R$-matrix.
For the case of the algebra $\zh(m_1,m_2,m_3,m_4)$, we start with the following $R$-matrix
\begin{equation}\label{eq:R}
 R(u)=\begin{pmatrix}
     uq -\frac{1}{uq} & 0 & 0 & 0 \\[2mm]
     0 &u-\frac{1}{u} &  q-\frac{1}{q}& 0\\[2mm]
     0 &q-\frac{1}{q} & u-\frac{1}{u}& 0\\[2mm]
     0&0&0&  uq -\frac{1}{uq}
      \end{pmatrix}.
\end{equation}
This $R$-matrix is associated to the quantum affine algebra
$U_q(\widehat{\mathfrak{sl}}_2)$ and is a solution of the Yang--Baxter equation
\begin{equation}
 R_{12}(u_1/u_2)R_{13}(u_1/u_3)R_{23}(u_2/u_3)
 =R_{23}(u_2/u_3)R_{13}(u_1/u_3)R_{12}(u_1/u_2),
\end{equation}
where $R_{12}=\sum_a R_a\otimes R^a\otimes\mathbb{1}$, $R_{23}=\sum_a\mathbb{1}\otimes
R_a\otimes R^a$, $R_{13}=\sum_a   R_a\otimes\mathbb{1}\otimes R^a$ if one writes $R=\sum_a
R_a\otimes R^a$ and $\mathbb{1}$ as the $2\times2$ identity matrix.  We define also the
following truncated reflection matrix (see remark \ref{rmk:trunc} below) given by
\begin{equation}\label{eq:K}
B(u)=\begin{pmatrix}
 \displaystyle  uqC_{12}-\frac{C_{23}}{uq}+\frac{p_4/u + p'_4u}{u^2-1/u^2}
 &  \displaystyle  qu^2 +\frac{1}{qu^2}-\frac{ [C_{23},C_{12}]_q}{q^2-1/q^2}
  +\frac{p''_4}{q+1/q}  \\[4mm]
 \displaystyle  - qu^2 -\frac{1}{qu^2}+\frac{ [C_{12},C_{23}]_q}{q^2-1/q^2}
  -\frac{p''_4}{q+1/q}
 & \displaystyle   uqC_{23}-\frac{C_{12}}{uq}+\frac{p_4u + p'_4/u}{u^2-1/u^2}
 \end{pmatrix},
\end{equation}
where we refer to \eqref{eq:p41}--\eqref{eq:p43} for the definition of $p_4$, $p'_4$ and
$p''_4$.

\begin{prop}\cite{Baseilhac2005}
The set of relations obtained from the reflection equation
\begin{equation}\label{eq:RKRK}
 R(u/v) B_1(u) R(uv) B_2(v)=B_2(v) R(uv) B_1(u) R(u/v),
\end{equation}
where $B_1(u)=B(u)\otimes\mathbb{1}$ and $B_2(u)=\mathbb{1}\otimes B(u)$, is equivalent to
the defining relations of $\zh(m_1,m_2,m_3,m_4)$.
\end{prop}
\proof We look at each matrix element of the reflection equation \eqref{eq:RKRK} and
derive $16$ relations. For each or them, we extract the different coefficients w.r.t. the
parameter $u$; this provides relations between $C_{12}$ and $C_{23}$. By direct
investigation, we verify that all the obtained relations are equivalent to the defining
relations of $\zh(m_1,m_2,m_3,m_4)$.
\endproof

Rephrasing this proposition, the Zhedanov algebra $\zh(m_1,m_2,m_3,m_4)$ is isomorphic to
the truncated reflection algebra defined by the $R$-matrix \eqref{eq:R} and the truncated
reflection matrix \eqref{eq:K}.
\begin{rmk}\label{rmk:trunc}
There exists a more general form for the reflection matrix, containing an infinite number
of generators encompassed in formal series of $u$ and $\frac1u$. The elements of the
reflection matrix \eqref{eq:K} can be obtained as a truncation of these formal series.
The algebra defined by the general reflection matrix obeying the reflection equation
\eqref{eq:RKRK} is isomorphic to the $q$-Onsager algebra \cite{BaseilhacKoizumi2005}.
Therefore, the Zhedanov algebra can also be seen as a quotient of the $q$-Onsager algebra.
\end{rmk}
In the context of the reflection algebra it is well-known how to obtain central elements
\cite{Sklyanin1988}. Indeed, let us define the Sklyanin determinant $\text{sdet} B(u)$ as
follows
\begin{equation}
 \mathrm{sdet} B(u):=-\frac12 tr_{12}\left(  R(1/q) B_1(u/q) R(u^2/q) B_2(u) \right).
\end{equation}
We can show that the coefficients of $\text{sdet} B(u)$ commute with $C_{12}$ and
$C_{23}$. We recover in this way that the operator $\varOmega$ given by \eqref{eq:aw3cas}
commutes with $C_{12}$ and $C_{23}$. The Sklyanin determinant gives solely $\varOmega$ as
a central element.  Fixing the Sklyanin determinant to an appropriate value allows us to
give a FRT presentation of $\szh(m_1,m_2,m_3,m_4)$:
\begin{prop}
The truncated reflection algebra defined by the $R$-matrix \eqref{eq:R}, the truncated
reflection matrix \eqref{eq:K} and quotiented by the relation
\begin{align}\label{eq:sd}
\begin{aligned}
 \mathrm{sdet} B(u)=q^2(1-q^4)^2&\,(u^2+q^{-m_2-m_4})
 (u^2+q^{m_2+m_4})(u^2+q^{m_4-m_2}) (u^2+q^{m_2-m_4})\\
 \times&\,(u^2+q^{-m_1-m_3}) (u^2+q^{m_1+m_3})(u^2+q^{m_3-m_1})
 (u^2+q^{m_1-m_3}),
\end{aligned}
\end{align}
is isomorphic to $\szh(m_1,m_2,m_3,m_4)$.
\end{prop}
\proof By direct computations, we show that \eqref{eq:sd} is equivalent to imposing
\eqref{eq:awt3cas1}. \endproof
The fact that $\szh(m_1,m_2,m_3,m_4)$ can be defined as a truncated reflection algebra was
expected, but it is a surprise that the r.h.s. of \eqref{eq:sd} factorizes into such a
simple form.

\section{A $\mathbf{W(D_4)}$ symmetry}\label{sec:szh}

The algebra $\szh(m_1,m_2,m_3,m_4)$ has a remarkable symmetry based on the
Weyl group $W(D_4)$ associated to the Lie algebra $D_4$. To describe it, let us
introduce a root system of type $D_4$ and fix a set of simple roots $\alpha_1$,
$\alpha_2$, $\alpha_3$, $\alpha_4$ with labeling according to the following Dynkin
diagram:
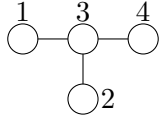
\begin{center}
\begin{tikzpicture}[scale=0.2]
\draw (5,1) circle [radius=1];
\node [above] at (5,1.5) {$1$};
\draw (6,1)--(8,1);

\draw (9,1) circle [radius=1];
\node [above] at (9,1.5) {$3$};
\draw (10,1)--(12,1);

\draw (13,1) circle [radius=1];
\node [above] at (13,1.5) {$4$};

\draw (9,-3) circle [radius=1];
\node [right] at (9.5,-3) {$2$};
\draw (9,0)--(9,-2);
\end{tikzpicture}
\end{center}
The Weyl group $W(D_4)$ is generated by the reflections $s_i$ associated to the simple
roots $\alpha_i$ which satisfy, for $1\leq i,j\leq 4$,
\begin{align}
\begin{aligned}
 s_i^2=&1\\
 s_is_j=&s_js_i\\
 s_is_js_i=&s_js_is_j
\end{aligned}\qquad
\begin{aligned}
 &\\
 &\text{if $i$ and $j$ are not connected in the Dynkin diagram},\\
 &\text{if $i$ and $j$ are connected in the Dynkin diagram}.
\end{aligned}
\end{align}
Its order is $192$.
Let us now associate the parameters $m_1,m_2,m_3,m_4$ with some of the roots as follows:
\begin{align}\label{}
 m_1=\alpha_1,\qquad m_2=\alpha_2,\qquad m_3=\alpha_4,\qquad m_4=\Theta,
\end{align}
where $\Theta$ is the longest positive root. The explicit expression of $\Theta$ is:
\begin{align}\label{}
 \alpha_3=\frac{1}{2}(m_4-m_1-m_2-m_3).
\end{align}
It is elementary to calculate the actions $s_i$ expressed in terms of the parameters:
\begin{equation}\label{eq:s}
s_1\,:\,m_1\mapsto -m_1,\quad\
s_2\,:\,m_2\mapsto -m_2,\quad\
s_4\,:\,m_3\mapsto -m_3,\quad\
s_3\,:\left\{\begin{array}{l}
m_1\mapsto m_1+\alpha_3\,,\\
m_2\mapsto m_2+\alpha_3\,,\\
m_3\mapsto m_3+\alpha_3\,,\\
m_4\mapsto m_4-\alpha_3\,,
\end{array}\right.\end{equation}
where the omitted actions are trivial and the explicit expression of $\alpha_3$ is
given above.  The action of the Weyl group is extended to any function as follows:
\begin{align}\label{}
 (\sigma f)(m_1,m_2,m_3,m_4)=f(\sigma(m_1),\sigma(m_2),\sigma(m_3),\sigma(m_4))
\end{align}
for $\sigma\in W(D_4)$.
\begin{prop}
The Weyl group $W(D_4)$ is a symmetry of $\szh(m_1,m_2,m_3,m_4)$ \textit{i.e.}
\begin{equation}
 \szh(m_1,m_2,m_3,m_4)=\szh(\sigma(m_1),\sigma(m_2),\sigma(m_3),\sigma(m_4)),
\end{equation}
for any $\sigma\in W(D_4)$.
\end{prop}
\proof In $\szh(m_1,m_2,m_3,m_4)$, we remark that the only functions of $m_i$
which appear are
\begin{subequations}
\begin{eqnarray}
 p_4&=& \chi_{m_1}\chi_{m_2}+\chi_{m_3}\chi_{m_4},\label{eq:p41}\\
 p'_4&=& \chi_{m_2}\chi_{m_3}+\chi_{m_1}\chi_{m_4},\\
 p''_4&=& \chi_{m_1}\chi_{m_3}+\chi_{m_2}\chi_{m_4},\label{eq:p43}\\
 p_6&=& \chi_{m_1}^2+\chi_{m_2}^2+\chi_{m_3}^2+\chi_{m_4}^2
  +\chi_{m_1}\chi_{m_2}\chi_{m_3}\chi_{m_4},
\end{eqnarray}
\end{subequations}
where $\chi_{m}=q^m+q^{-m}$.
By direct computations, we can show that these functions are invariant by the
transformations $s_1$, $s_2$, $s_3$ and $s_4$ given by \eqref{eq:s}, which concludes the
proof since they generate $W(D_4)$.
\endproof
In the study of the finite representations of the universal algebra $\Delta_q$ a $W(D_4)$
symmetry has been also investigated \cite{Terwilliger2004, NomuraTerwilliger2009}.

\subsection{Connection with the $W(D_4)$ symmetry in the Racah algebra}

Let us perform the transformation
\begin{align}
 K_{I}=\frac{C_I-(q+q^{-1})}{(q-q^{-1})^{2}},
\end{align}
with $I \in \{ 1,2,3,123,12,23\}$. Note that $13$ does not belong to this set.
In the algebra $\szh(m_1,m_2,m_3,m_4)$, one gets, for $i=1,2,3,4$,
\begin{equation}
 K_i = \frac{\chi_{m_i}-(q+q^{-1})}{(q-q^{-1})^{2}}
 =\left[\frac{m_i}{2}\right]_q^2-\left[\frac{1}{2}\right]_q^2,
\end{equation}
where the $q$-number is defined by $[m]_q=\frac{q^m-q^{-m}}{q-q^{-1}}$.
The commutation relations of the algebra $\szh(m_1,m_2,m_3,m_4)$ become
\begin{subequations}
\begin{align}
 & [K_{12},K_{23}]_q=K_{13}, \label{eq:K1}\\
 & [K_{23},K_{13}]_q=(q+q^{-1})\big(- \{K_{12}, K_{23} \}  - K_{23}^2
  + \xi_2 K_{23} +\xi_4\big),\label{eq:K2}\\
 & [K_{13},K_{12}]_q=(q+q^{-1})\big(- \{K_{12}, K_{23} \}  - K_{12}^2
  + \xi_2 K_{12} +\xi'_4\big) ,\label{eq:K3}
\end{align}
and the supplementary relation becomes
\begin{equation}
\begin{aligned}
&-q\frac{q-q^{-1}}{q+q^{-1}}\,K_{12}K_{23}K_{13}-qK_{12}K_{23}K_{12}
-q^{-1}K_{23}K_{12}K_{23}
+\frac{q^2}{(q+q^{-1})^2}\,K_{13}^2\\
&+\Big( \frac{\xi_2}{q+q^{-1}}-1\Big)\{K_{12}, K_{23}\}+q\xi_4
K_{12}+q^{-1}\xi'_4 K_{23}
\quad=\quad
\xi_6 - \xi_4 -\xi'_4-\frac{\xi_2^2}{4},
\end{aligned}
\end{equation}
with
\begin{align}
 & \xi_2= \frac{1}{q+q^{-1}}\Big( 2(M_1^2 + M_2^2 + M_3^2+ M_4^2 - 1)
 + (q-q^{-1})^2(M_1^2 M_3^2+M_2^2M_4^2)\Big),\label{eq:xxi2}\\
 & \xi_4= (M_1^2-M_{4}^2)(M_3^2-M_2^2),\\
 & \xi_4'= (M_1^2-M_{2}^2)(M_{3}^2-M_4^2),\label{eq:xxi4p}\\
 & \xi_6= ( M_1^2 M_3^2 - M_2^2M_4^2)(M_1^2-M_2^2+M_3^2 -M_4^2)
  + \tfrac{1}{4}(q-q^{-1})^2(M_1^2 M_3^2 - M_2^2M_4^2)^2,\label{eq:xxi6}
\end{align}
\end{subequations}
where we use the notation $M_i=\left[\frac{m_i}{2}\right]_q $.  As expected, we can check
that the functions $\xi_2$, $\xi_4$, $\xi'_4$ and $\xi_6$ are invariant under the action
of the Weyl group $W(D_4)$.

The advantage of this presentation of $\szh(m_1,m_2,m_3,m_4)$ is that the classical limit
$q\to 1$ (see Appendix \ref{app:inj}) is well-defined and provides straightforwardly the
commutation relations of the Racah algebra.  Thus, the description of the Weyl group
$W(D_4)$ action also holds for the Racah algebra and we recover the results of
\cite{GranovskiiZhedanov1988}.
\begin{rem}
In the classical limit $q\to1$, the functions $\xi_2$, $\xi_4$, $\xi'_4$ and $\xi_6$ form
a basis for polynomials invariant under the action of $W(D_4)$, as expected. In the
generic case ($q\in\mathbb{C}$, not a root of unity), one gets two different sets of
invariant functions: $S_{\xi}=\{\xi_2, \xi_4, \xi'_4, \xi_6\}$ on the one hand and
$S_p=\{p_4,p'_4,p''_4,p_6\}$ on the other hand. We have checked that there exists an
invertible polynomial mapping between these two sets.  However, only $S_{\xi}$ admits a
non-trivial classical limit.
\end{rem}

\section{Kauffman bracket skein modules and algebras}\label{sec:KBS}

Kauffman bracket skein module quantizations have been introduced in
\cite{BullockPrzytycki1999, Turaev1991} and further studied along our lines of interest
for this paper in \cite{BerestSamuelson2018, Hikami2019, CookeLacabanne2020}.  We will now
recall some key definitions and results from these investigations. We shall work with an
oriented $3$-manifold $\mathcal{M}$ which is a thickened surface, that is
$\mathcal{M}=\Sigma_{0,n}\times I$, where $I=[0,1]$ and $\Sigma_{0,n}$ is the
$n$-punctured sphere.

\begin{defi}
The quantized skein module $\Sk_{\theta}(\mathcal{M})$ is the
$\mathbb{C}[\theta^{\pm1}]$-module spanned by framed and unoriented links in $\mathcal{M}$
modulo the Kauffman bracket skein relations that allow to ``simplify the crossings'':
\begin{subequations}\label{eq:Kauf}
\begin{align}
\KPB&= \theta \KPC+\theta^{-1}\KPD,\label{eq:KaufA}\\[.5em]
\KPA&= -(\theta^{2}+\theta^{-2}),\label{eq:KaufB}
\end{align}
\end{subequations}
where $\theta\in\mathbb{C}$ is not a root of unity and in the framing relation
\eqref{eq:KaufB} the link should not enclose a puncture.  This defines an algebra, which
we will denote $\Sk_{\theta}(\Sigma_{0,n})$, for which multiplication is given by stacking
the links on top of each other in the $I$ direction.
\end{defi}\noi
We shall use diagrams that correspond to the projection of the links on the surface (all
the while keeping the information about the relative ``height'' of the links in the $I$
direction). Let us now establish the conventions for these drawings (framed links
diagrams).

The $n$-punctured sphere $\Sigma_{0,n}$ is equivalent to the plane with $n-1$ punctures
(denoted by the $(n-1)$ drawn $\times$'s):
\begin{align}\label{fig:dashcont}
&\begin{tikzpicture}
\draw (0,0) node{$\times$};
\draw (1,0) node{$\times$};
\draw (2,0) node{$\times$};
\draw (3,0) node{$\dots$};
\draw (4,0) node{$\times$};
\draw[dashed] (-1,-1) -- (5,-1);
\draw[dashed] (-1,1) -- (5,1);
\draw[dashed] (-1,-1) -- (-1,1);
\draw[dashed] (5,-1) -- (5,1);
\end{tikzpicture}
\end{align}
The dashed contour corresponds to the $n$\textsuperscript{th} puncture of the sphere. We
will omit the contour in the subsequent diagrams but it is always understood to be there.

Framed links that enclose punctures are represented by loops drawn around the $\times$'s.
We shall use the term ``loops'' to refer unambiguously to the framed links in the
remainder of the paper. These loops can be homotopically deformed without crossing the
holes (punctures). Remark that loops enclosing a single puncture are central elements in
$\Sk_{\theta}(\Sigma_{0,n})$. This is also true for the $n^{th}$ puncture, which amounts
to saying that the loop enclosing the $(n-1)$ punctures $\times$ is also central.

Let us now consider the surface $\Sigma_{0,4}$ and give names to a few loops:
\begin{align}\label{eq:curvesAW3}
\begin{aligned}
&\begin{tikzpicture}
\punctiii \larc{0} \rarc{1} \updown{0}{1} \draw (2.5,0)
node[label=right:{$=\mathbb{A}_{12}$}]{};
\end{tikzpicture}\\[.em]
&\begin{tikzpicture}
\punctiii \larc{1} \rarc{2} \updown{1}{2} \draw (2.5,0)
node[label=right:{$=\mathbb{A}_{23}$}]{};
\end{tikzpicture}\\[.em]
&\begin{tikzpicture}
\punctiii \uarc{2} \uarc{0} \bendd{0.25}{1.75} \bendd{-0.25}{2.25} \draw (2.5,0)
node[label=right:{$=\mathbb{A}_{13}$}]{};
\end{tikzpicture}
\end{aligned}\hspace{6em}
\begin{aligned}
&\begin{tikzpicture}
\punctiii \larc{0} \rarc{0} \draw (2.5,0) node[label=right:{$=\mathbb{A}_{1}$}]{};
\end{tikzpicture}\\[.em]
&\begin{tikzpicture}
\punctiii \larc{1} \rarc{1} \draw (2.5,0) node[label=right:{$=\mathbb{A}_{2}$}]{};
\end{tikzpicture}\\[.em]
&\begin{tikzpicture}
\punctiii \larc{2} \rarc{2} \draw (2.5,0) node[label=right:{$=\mathbb{A}_3$}]{};
\end{tikzpicture}\\[.em]
&\begin{tikzpicture}
\punctiii \larc{0} \rarc{2} \updown{0}{2} \draw (2.5,0)
node[label=right:{$=\mathbb{A}_{123}$}]{};
\end{tikzpicture}
\end{aligned}
\end{align}
Following the definition, multiplication of two loops $X\cdot Y$ means putting $Y$
\emph{on top of} $X$, for example:
\begin{align}
\begin{aligned}
&\begin{tikzpicture}
\punctiii \larc{0} \larc{1} \rarc{2} \updown{1}{2}
\draw (0,-0.25cm) -- (0.8,-0.25cm); \draw (0,0.25cm) -- (0.8,0.25cm);
\draw (1.cm,-0.2cm) arc (-80:80:0.2cm);
\draw (-0.5,0) node[label=left:{$\mathbb{A}_{12}\cdot\mathbb{A}_{23}=$}]{};
\end{tikzpicture}\\[.em]
\end{aligned}
\end{align}
One would then proceed to use relations \eqref{eq:Kauf} to simplify the expressions:
\begin{align}
\begin{aligned}
 \mathbb{A}_{12}\cdot \mathbb{A}_{23}=&
 ~\theta
 \left(\begin{aligned}\begin{tikzpicture}
 \punctiii
 \draw (0.5cm,-0.25cm) -- (0.cm,-0.25cm) arc (270:90:0.25cm) -- (0.8cm,0.25cm)
 arc (90:-90:0.05cm) arc (90:270:0.2cm) -- (2cm,-0.25cm)
 arc (-90:90:0.25cm) -- (1.1cm,0.25cm) arc (90:270:0.05cm) arc (90:-90:0.15cm)
 -- (0.75cm,-0.15cm) arc (90:125:0.1cm);
 \draw (0.5cm,-0.25cm) arc (-90:-30:0.1cm);
 \end{tikzpicture}\end{aligned}\right)
 -\theta^{-1}
 \left(\begin{aligned}\begin{tikzpicture}
 \punctiii
 \draw (0.7cm,-0.25cm) -- (0.cm,-0.25cm) arc (270:90:0.25cm) -- (2cm,0.25cm)
 arc (90:-90:0.25cm) -- (1cm,-0.25cm) arc (270:90:0.225cm)
 arc (90:-90:0.2cm) -- (0.95cm,-0.2cm);
 \draw (0.7cm,-0.25cm) arc (-90:-30:0.1cm);
 \end{tikzpicture}\end{aligned}\right)\\[2em]
 =&\phantom{+}\theta^{2}
 \left(\begin{aligned}\begin{tikzpicture}
 \punctiii \uarc{2} \uarc{0} \bendd{0.25}{1.75} \bendd{-0.25}{2.25}
 \end{tikzpicture}\end{aligned}\right)\qquad
 +\phantom{\theta^{-2}}\left(\begin{aligned}\begin{tikzpicture}
 \punctiii \larc{0}\rarc{2}\updown{0}{2} \draw (1cm,0.2cm) arc (90:450:0.2cm);
 \end{tikzpicture}\end{aligned}\right)\\[0em]
 &+\phantom{\theta^{2}}
 \left(\begin{aligned}\begin{tikzpicture}
 \punctiii \larc{0}\rarc{0} \larc{2}\rarc{2}
 \end{tikzpicture}\end{aligned}\right)\qquad
 +\theta^{-2}
 \left(\begin{aligned}\begin{tikzpicture}
 \punctiii \darc{2} \darc{0} \bendu{0.25}{1.75} \bendu{-0.25}{2.25}
 \end{tikzpicture}\end{aligned}\right).\\
  =& ~\theta^{2}\mathbb{A}_{13} + \mathbb{A}_2\cdot\mathbb{A}_{123}
   + \mathbb{A}_1\cdot\mathbb{A}_{3}
   +\theta^{-2}
 \left(\begin{aligned}\begin{tikzpicture}
 \punctiii \darc{2} \darc{0} \bendu{0.25}{1.75} \bendu{-0.25}{2.25}
 \end{tikzpicture}\end{aligned}\right).
\end{aligned}
\end{align}
Similarly, exchanging the order of multiplication, one obtains the same diagrams but with
inverse coefficients:
\begin{align}
\begin{aligned}
 \mathbb{A}_{23}\cdot \mathbb{A}_{12}
 =& ~\theta^{-2} \mathbb{A}_{13} +\mathbb{A}_2\cdot\mathbb{A}_{123}
   + \mathbb{A}_1\cdot\mathbb{A}_{3}
   +\theta^{2}
 \left(\begin{aligned}\begin{tikzpicture}
 \punctiii \darc{2} \darc{0} \bendu{0.25}{1.75} \bendu{-0.25}{2.25}
 \end{tikzpicture}\end{aligned}\right).
\end{aligned}
\end{align}
We see immediately that one gets
\begin{equation}
\theta^2\mathbb{A}_{12}\cdot\mathbb{A}_{23}-\theta^{-2}\mathbb{A}_{23}\cdot\mathbb{A}_{12}
 =(\theta^4-\theta^{-4}) \mathbb{A}_{13}+
 (\theta^2-\theta^{-2})(\mathbb{A}_2\cdot \mathbb{A}_{123}
 +\mathbb{A}_1\cdot \mathbb{A}_{3}).
\end{equation}
The skein algebra $Sk_{\theta}(\Sigma_{0,4})$ is directly linked to the Askey--Wilson
algebra as stated in the following proposition:
\begin{prop}\label{pr:as}
The Special Askey--Wilson algebra $\saw(3)$ is isomorphic to the Kauffman bracket skein
algebra $\Sk_{\,iq^{1/2}}(\Sigma_{0,4})$.  The isomorphism is given by the following
invertible map:
\begin{align}\label{eq:equivaw3tS04}
C_{I} &\mapsto \mathbb{A}_{I},
\end{align}
for $I \in \{1,2,3,123,12,23,13\}$.
\end{prop}
\proof The isomorphism is directly verified by comparing the relations of $\saw(3)$ and
the ones of the Kauffman bracket skein algebra obtained in \cite{BullockPrzytycki1999}
(see also Proposition $3.1$ of \cite{Hikami2019} and \cite{Cooke2020} for additional
details). \endproof

This proposition gives a diagrammatic approach to study the algebra $\saw(3)$.\\[1em]
Let us emphasize that the previous isomorphism involves the Special Askey--Wilson algebra
$\saw(3)$. If we replace $\saw(3)$ by $\mathbf{aw}(3)$ in the map of the proposition, the
homomorphism would be not injective and if we instead replace $\saw(3)$ by $\Delta_q$ (as
in \cite{Terwilliger2015, Huang2016}), it would be not surjective.

One notes that the $\mathbb{Z}_3$-symmetry of the $\saw(3)$ relations is made manifest in
terms of the framed links picture, as the punctures do not have fixed positions and can be
switched around.

From now on we will unambiguously refer to the drawn loops identified as the generators of
$Sk_{\,iq^{1/2}}(\Sigma_{0,4})$ directly as their $C_{I}$ counterpart following
\eqref{eq:equivaw3tS04}.  This correspondence \eqref{eq:equivaw3tS04} leads to a natural
labeling of the punctures.  Indeed, consider the generators given in
\eqref{eq:curvesAW3}: the punctures enclosed in a given loop correspond precisely to the
set of indices $I$ of the corresponding generator $C_{I}$ if one labels the punctures
consecutively as:
\begin{align}\label{eq:label3punct}
\begin{tikzpicture}
\punctiii
\draw (0,0) node[label=below:{\scriptsize$1$}]{};
\draw (1,0) node[label=below:{\scriptsize$2$}]{};
\draw (2,0) node[label=below:{\scriptsize$3$}]{};
\end{tikzpicture}
\end{align}
\begin{rmk}
We recall that one arrives to the Special Zhedanov algebra $\szh(m_1,m_2,m_3,m_4)$ from
the Special Askey--Wilson algebra $\saw(3)$ by attributing a value to the central elements
$C_i$, $i=1,2,3,123$, see \eqref{eq:mm}. In the same way, starting from the Kauffman
bracket skein algebra $\Sk_{\,iq^{1/2}}(\Sigma_{0,4})$, one can define an evaluated
Kauffman bracket skein algebra, denoted ${\Sk}_{\,iq^{1/2}}(\Sigma_{0,4};m_1,m_2,m_3,m_4)$
by attributing a value to the puncture-framing relations:
\begin{align}\label{eq:skein-arc}
\begin{aligned}
&\begin{tikzpicture}
\draw (0,0) node{\phantom{$\times$}};
\draw (1,0) node{$\times$}; \larc{1} \rarc{1};
\draw (1,0) node[label=below:{\scriptsize{$i$}}]{};
\draw (2.5,0) node[label=right:{$=q^{m_i}+q^{-m_i}\,,\quad i=1,2,3,$}]{};
\end{tikzpicture}
\\
&\begin{tikzpicture}
\punctiii \larc{0} \rarc{2} \updown{0}{2}
\draw (2.5,0) node[label=right:{$=q^{m_4}+q^{-m_4}$.}]{};
\end{tikzpicture}
\end{aligned}
\end{align}
Note that the last drawing corresponds in fact to a contour enclosing the fourth puncture
on the sphere, see \eqref{fig:dashcont}.  As a corollary of Proposition \ref{pr:as}, the
algebra ${\Sk}_{\,iq^{1/2}}(\Sigma_{0,4}; m_1,m_2,m_3,m_4)$ is isomorphic to the Special
Zhedanov algebra $\szh(m_1,m_2,m_3,m_4)$.

Relations \eqref{eq:skein-arc} with $m_i=1$ already appear in the definition of the skein
algebra of arcs and link introduced in \cite{RogerYang2014}, from where we borrowed the
terminology `puncture-framing'.
\end{rmk}

\section{$U_q(\sl_2)$ and its centralizer in $U_q(\sl_2)^{\otimes3}$}\label{sec:cent}

The goal of this section is to discuss the notion of centralizer of $U_q(\sl_2)$ in
$U_q(\sl_2)^{\otimes3}$, which we denote by $\mC_3$, and connect it with the Special
Askey--Wilson algebra $\saw(3)$.

\subsection{$U_q(\sl_2)$ and its universal $R$-matrix\label{sect:Uqsl2}}

Let us fix the notation and conventions that will be used to perform the explicit
calculations in $U_q(\sl_2)$ (note that the results obtained will be independent of these
conventions at the end). We shall first define the quasi-triangular Hopf algebra
$U_q(\sl_2)$, present its braided universal $R$-matrix and list some additional properties
of interest.

$U_q(\sl_2)$ is an associative algebra generated by $E$, $F$,  $q^{H}$ and $q^{-H}$
obeying the defining relations
\begin{align}
 q^{H}q^{-H}=q^{-H}q^{H}=1, \quad q^{H}E=q Eq^H, \quad q^{H}F=q^{-1} Fq^H  \quad
 \text{and} \quad [E,F]=[2H]_q.
\end{align}
The center of this algebra is generated by the following Casimir element
(denoted $\Lambda$ in \cite{DeBieDeClercqetal2020, DeClercq2019})
\begin{align}\label{eq:Qcas}
Q=(q-q^{-1})^2\left(FE + \frac{qq^{2H}+q^{-1}q^{-2H}}{(q-q^{-1})^2}\right).
\end{align}
The algebra $U_q(\sl_2)$ can be endowed with a Hopf structure. In
particular, its comultiplication (or coproduct) homomorphism $\Delta:U_q(\sl_2)
\rightarrow U_q(\sl_2)\otimes U_q(\sl_2) $ is given by
\begin{subequations}
\begin{align}
 \Delta(E)&=E \otimes q^{-H} +q^H \otimes E,\qquad
 \Delta(q^H)=q^H\otimes q^H,\\
 \Delta(F)&=F \otimes q^{-H} +q^{H} \otimes F,\qquad
 \Delta(q^{-H})=q^{-H}\otimes q^{-H},
\end{align}
\end{subequations}
and is coassociative
\begin{align}\label{eq:com}
 (\Delta \otimes \id)\Delta= (\id \otimes \Delta)\Delta.
\end{align}
The quantum algebra $U_q(\sl_2)$ is called quasi-triangular because in a completion of
$U_q(\sl_2)\otimes U_q(\sl_2)$, there exists a universal $R$-matrix $\cR$ which is
invertible and satisfies
\begin{align}
 \Delta(x) \cR  &= \cR \Delta^{op}(x) \qquad \text{for } x\in U_q(\sl_2),\label{eq:RD} \\
 (\id \otimes \Delta)\cR&=\cR_{12}\cR_{13},\\
  (\Delta \otimes \id)\cR&=\cR_{23}\cR_{13},
\end{align}
where in the Sweedler notation we write the opposite comultiplication
$\Delta^{op}(x)=x_{(2)} \otimes x_{(1)}$ if $\Delta(x)=x_{(1)} \otimes x_{(2)}$.
In the previous relation, we have used the notations
$\cR_{12}=\cR^\alpha\otimes \cR_\alpha \otimes 1$, $\cR_{23}=1\otimes \cR^\alpha\otimes
\cR_\alpha$ and $\cR_{13}=\cR^\alpha\otimes1\otimes\cR_\alpha$ where
$\cR=\cR^\alpha\otimes \cR_\alpha$ (the sum over repeated indices $\alpha$ is understood).
The universal $R$-matrix is
given explicitly by \cite{Drinfeld1988}
\begin{equation}\label{eq:uR}
 \cR=q^{2(H\otimes H)} \sum_{n=0}^\infty \frac{(q-q^{-1})^n}{[n]_q!} q^{n(n-1)/2}
 \big(E q^H  \otimes q^{-H} F\big)^n,
\end{equation}
where $[n]_q!=[n]_q[n-1]_q\dots [2]_q[1]_q$ and, by convention, $[0]_q!=1$.

One can also define the so-called braided universal $R$-matrix $\check\cR$ by
\begin{align}\label{eq:czR}
 \czR_{i}=\cR_{i,i+1}\sigma_{i,i+1}
\end{align}
where $\sigma_{i,i+1}$ acts on the $i$\textsuperscript{th} and $(i+1)$\textsuperscript{th}
factors of the tensor product as
\begin{align}
 \sigma_{i,i+1}(\dots\otimes x_i\otimes x_{i+1}\otimes \dots)
 =(\dots \otimes x_{i+1}\otimes x_{i}\otimes\dots)\sigma_{i,i+1}.
\end{align}
This braided universal $R$-matrix satisfies the braided Yang--Baxter equation
\begin{equation}\label{bYBE}
 \czR_{i} \czR_{i+1} \czR_{i}=\czR_{i+1} \czR_{i} \czR_{i+1}.
\end{equation}

\subsection{An algebra generated by the intermediate Casimir elements}

Let us define the following \textit{intermediate} Casimir elements
\begin{gather}\label{eq:Casim3}
\begin{gathered}
 Q_1=Q\otimes 1\otimes 1,\qquad Q_2=1\otimes Q\otimes 1,\qquad Q_3=1\otimes 1\otimes Q,\\
 Q_{12}=\Delta(Q) \otimes 1=Q_{(1)} \otimes Q_{(2)} \otimes 1,
 \qquad Q_{23}=1 \otimes \Delta(Q) =1 \otimes Q_{(1)} \otimes Q_{(2)},\\
 Q_{123}=(\Delta\otimes\id)\Delta(Q).
\end{gathered}
\end{gather}
The labeling of these intermediate Casimir elements is chosen so as to refer to the
non-trivial factors in the tensor product $U_q(\sl_2)^{\otimes 3}$.
\begin{defi}
The algebra $\mathcal{A}_3$ is the subalgebra of $U_q(\sl_2)^{\otimes 3}$ generated by the
intermediate Casimir elements $Q_1$, $Q_2$, $Q_3$, $Q_{12}$, $Q_{23}$ and $Q_{123}$.
\end{defi}
Let us define an additional intermediate Casimir element
\begin{equation}
Q_{13}=\czR_{2}^{-1}Q_{12}\czR_{2}=\czR_{1}Q_{23}\czR_{1}^{-1}.
\end{equation}
It has been proven in \cite{CrampeGaboriaudetal2020} that this element is in
$U_q(\sl_2)^{\otimes 3}$ (and not in its completion), that the second equality is
compatible with the first one and that the following proposition holds:
\begin{prop}
The intermediate Casimir elements $Q_1$, $Q_2$, $Q_3$, $Q_{123}$, $Q_{12}$, $Q_{23}$ and
$Q_{13}$ belong to the centralizer $\mC_3$ of the diagonal action of $U_q(\sl_2)$ in
$U_q(\sl_2)^{\otimes 3}$ defined by
\begin{equation}\label{eq:defC3}
 \mC_3=\{ X \in U_q(\sl_2)^{\otimes 3}\ \big| \ [(\Delta\otimes \id)\Delta(x), X]=0\ ,
 \ \  \forall x\in U_q(\sl_2) \} .
 \end{equation}
\end{prop}
The precise links between the Askey--Wilson algebra, the centralizer and the algebra
$\mathcal{A}_3$ generated by the intermediate Casimir elements are given in the following
proposition.
\begin{prop}\label{prop:CQ}
The algebra $\saw(3)$ has an homomorphic injective image in $\mC_3$. The mapping is done
as follows:
\begin{align}\label{eq:cQQ}
 C_I\mapsto Q_I,\qquad\text{for}\qquad I\in\{1,2,3,123,12,23,13\}.
\end{align}
The algebra $\saw(3)$ is isomorphic to $\mathcal{A}_3$.
\end{prop}
\proof All the relations of $\saw(3)$ given by \eqref{eq:aw3} and \eqref{eq:awt3cas1} are
easily checked in $U_q(\sl_2)^{\otimes 3}$ upon rewriting the $Q_I$'s in terms of the
$U_q(\sl_2)^{\otimes 3}$ generators. The proof of the injectivity is postponed to Appendix
\ref{app:inj}. {The method used in \cite{Huang2016} to prove the injectivity of
$\Delta_q$ into $U_q(\sl_2)^{\otimes 3}$ seems difficult to generalize to the case treated
here and we propose an alternative method based on classical invariant theory.} Since the
algebra $\mathcal{A}_3$ is the image of the map \eqref{eq:cQQ}, it follows that $\saw(3)$
is isomorphic to $\mathcal{A}_3$. \endproof

This realization of the Askey--Wilson algebra in $U_q(\sl_2)^{\otimes 3}$ was the
motivation for adding the relation \eqref{eq:awt3cas1} to the ``intuitive'' set of
relations of $\mathbf{aw}(3)$. Indeed, since relation \eqref{eq:awt3cas1} is obeyed by the
intermediate Casimir elements, it should also be included in the algebra encoding the
properties of these Casimir elements.
\begin{coro}\label{co:CQ}
The algebra $\mathcal{A}_3$ is isomorphic to the Kauffman bracket skein algebra
of the four-punctured sphere
$\Sk_{\,iq^{1/2}}(\Sigma_{0,4})$. The isomorphism is given by the following map:
\begin{align}\label{eq:cSQ}
 \phi : Q_I \mapsto \mathbb{A}_I,\qquad\text{for}\quad I\in\{1,2,3,123,12,23,13\}.
\end{align}
\end{coro}
\proof A direct consequence of the Propositions \ref{pr:as} and \ref{prop:CQ}.

\subsection{Fundamental theorems of invariant theory}\label{sec:Zmm}

In the previous section, we introduced the centralizer $\mC_3$ of the diagonal action of
$U_q(\sl_2)$ in the threefold tensor product and showed its connection with the
Askey--Wilson algebra $\saw(3)$.  We now focus on similar objects in the case where
we represent each factor $U_q(\sl_2)$ in $U_q(\sl_2)^{\otimes 3}$ by a finite-dimensional
irreducible representation.

The quantum algebra $U_q(\sl_2)$ has finite irreducible representations of dimension
$m=2j+1$ that we will denote by $M(m)$, with $m\in \mathbb{Z}_{>0}$.  The name ``spin-$j$
representation'' is usually used to refer to $M(m=2j+1)$.  The representation map will be
denoted by $\pi_{m}:U_q(\sl_2)\rightarrow \text{End}(M(m))$.  The representation of the
Casimir element \eqref{eq:Qcas} in the space $M(m)$ is
\begin{equation}
\pi_{m}(Q)= \chi_m \mathbb{1}_{m},
\label{eq:casimirrep}
\end{equation}
where $\chi_m=q^{m}+q^{-m}$ and $\mathbb{1}_m$ is the $m\times m$ identity matrix.

From now on, we fix three integers $m_1$, $m_2$ and $m_3$.
The threefold tensor product of irreducible representations of $U_q(\sl_2)$ decomposes
into the following direct sum of irreducible representations
\begin{equation}\label{eq:decomp2}
M(m_1) \otimes M(m_2) \otimes M(m_3)=\bigoplus_{m_4} M(m_4) \otimes V_{m_1,m_2,m_3}^{m_4},
\end{equation}
where $V_{m_1,m_2,m_3}^{m_4}$ is called the multiplicity space.  We recall that we look at
cases where $q$ is not a root of unity otherwise the previous statement would be wrong.

We now fix four integers $m_1$, $m_2$, $m_3$, $m_4$ and denote by $\cQ_I$ the image of
$Q_I$ in $ V_{m_1,m_2,m_3}^{m_4}$ (for $I\in \{1,2,3,123,12,23,13\}$).  We get
$\cQ_1=\chi_{m_1}$, $\cQ_2=\chi_{m_2}$, $\cQ_3=\chi_{m_3}$ and $\cQ_{123}=\chi_{m_4}$.
\begin{prop} There exists a surjective algebra homomorphism from $\szh(m_1,m_2,m_3,m_4)$
to $End(V_{m_1,m_2,m_3}^{m_4})$ given by
\begin{equation}
 C_I \mapsto \cQ_I, \qquad \text{for } I\in\{12,23,13\}.
\end{equation}
\end{prop}
This proposition which provides the generators for the centralizer of the diagonal action
is sometimes called in invariant theory the ``first fundamental theorem''.  The map in the
previous proposition is not injective. The description of the kernel of this map is the
subject of \cite{CrampeVinetetal2020a} (see also \cite{FloresPeltola2020}) and is called
the ``second fundamental theorem''.

We recall that the algebra $\szh(m_1,m_2,m_3,m_4)$ possesses a $W(D_4)$-symmetry.  Let us
remark that a similar Weyl group symmetry of type $E_6$ has been discovered recently
\cite{CrampePoulaindAndecyetal2020} in the case of the centralizer of the diagonal
embedding of $U(\sl_3)$ in two copies of $U(\sl_3)$.

\section{The Double Affine Hecke Algebra $(C_1^{\vee},C_1)$}\label{sec:daha}

Double affine Hecke algebras (DAHA) of type $(C_1^{\vee},C_1)$ were introduced in
\cite{Sahi1999} and their connections with Askey--Wilson polynomials were first explored
in \cite{NoumiStokman2004} and \cite{Oblomkov2004}. Universal analogues of these DAHA were
later introduced and studied in \cite{ItoTerwilliger2010, Terwilliger2013, Huang2020}.

In this section, we present another connection between the Special Askey--Wilson algebra
$\saw(3)$ and a certain subalgebra of a universal DAHA of type $(C_1^{\vee},C_1)$.
\begin{defi}\label{def:UDAHAPT}
We introduce the following algebras
\begin{itemize}
\item The universal Double Affine Hecke Algebra of type $(C_1^{\vee},C_1)$
\cite{ItoTerwilliger2010} is defined as the associative algebra $\Hq$ with generators
$\{t_i^{\pm1}\,,\ i=0,\dots,3\}$ and relations:
\begin{subequations}\label{eq:UDAHA}
\begin{align}
 &t_i t_i^{-1}=t_i^{-1}t_i=1,\\
 &t_i+t_i^{-1}\qquad\text{is central},\\
 &t_0t_1t_2t_3=q^{-1}.
\end{align}
\end{subequations}
The ``usual'' DAHA, denoted $H_q(k_0,k_1,k_2,k_3)$, is recovered when the central elements
$t_i+t_i^{-1}$ have complex values $k_i+k_i^{-1}$, with $k_i\neq0$.
\item The algebra $\Gamma_q$ \cite{Terwilliger2013} is the subalgebra of $\Hq$ commuting
with the distinguished generator $t_0$ ($\Gamma_q$ is the centralizer of $t_0$ in $\Hq$):
\begin{align}\label{eq:defiAq}
 \Gamma_q=\{h\in\Hq~|~[h,t_0]=0\}.
\end{align}
\item Let $\mathbf{e}$ be the following idempotent of $H_q(k_0,k_1,k_2,k_3)$
\cite{Hikami2019}
\begin{align}\label{eq:idempotent}
 \mathbf{e}=\frac{t_0-k_0}{k_0^{-1}-k_{0}}.
\end{align}
The spherical DAHA, denoted $\SHq(k_0,k_1,k_2,k_3)$
\cite{Koornwinder2007a,Koornwinder2008}, is defined as
\begin{align}\label{}
 \SHq(k_0,k_1,k_2,k_3)=\mathbf{e}\,H_q(k_0,k_1,k_2,k_3)\,\mathbf{e}.
\end{align}
\end{itemize}
\end{defi}

The following theorems relate DAHA to the previously introduced algebraic structures.
\begin{thm}\label{thm:terwil}
\cite{Terwilliger2013} The map $\Theta:\ \saw(3)\to\Gamma_q$ defined by
\begin{align}
\begin{aligned}
C_{12}&\mapsto t_1t_0+(t_1t_0)^{-1} ,\\
C_{23}&\mapsto t_3t_0+(t_3t_0)^{-1} ,\\
C_{13}&\mapsto t_2t_0+(t_2t_0)^{-1},
\end{aligned}\qquad\quad
\begin{aligned}
C_1&\mapsto t_1+t_1^{-1},\\
C_2&\mapsto t_2+t_2^{-1},\\
C_3&\mapsto t_3+t_3^{-1},\\
C_{123}&\mapsto q^{-1}t_0 + qt_0^{-1}.
\end{aligned}
\end{align}
is an injective algebra homomorphism.
\end{thm}

\begin{thm}\label{thm:koorn} (Theorem 3.2 in \cite{Koornwinder2008})
The Special Zhedanov algebra $\szh(m_1,m_2,m_3,m_4)$ is isomorphic to the spherical DAHA
$\SHq(k_0,k_1,k_2,k_3)$.
\end{thm}
\begin{rmk}
Spherical DAHAs have also been connected to skein algebras of higher genus. The Kauffman
bracket skein algebra of the once-punctured torus $\Sk_\theta(\Sigma_{1,1})$ is related
to a (spherical) DAHA of type $A_1$ \cite{BullockPrzytycki1999, Samuelson2019} and the
genus two skein algebra is related to a genus two spherical double affine Hecke algebra in
\cite{CookeSamuelson2020}.
\end{rmk}

\section{Actions of the braid group}\label{sec:braid}

In this section, we provide two actions of the braid group: the first one on the algebra
$\mathcal{A}_3$ and the second one on the skein algebra $\Sk_{\,iq^{1/2}}(\Sigma_{0,4})$.
Then, we show how these two actions are compatible and give a diagrammatic presentation of
the intermediate Casimir elements of $U_q(\sl_2)^{\otimes 3}$.

We recall that the braid group on $n$ strands $B_n$ is generated by the elements
$s_1,\dots,s_{n-1}$ as well as their inverses $s_1^{-1},\dots,s_{n-1}^{-1}$ satisfying
\begin{align}\label{eq:braid_relations}
\begin{aligned}
 s_is_{i+1}s_i&=s_{i+1}s_is_{i+1},\\
 \quad s_is_j&=s_js_i\qquad\quad\text{if} \quad |i-j|\geq2,\\
 s_i^{-1}s_i&=s_is_i^{-1}=1.
\end{aligned}
\end{align}

\subsection{The braided universal $R$-matrix and a braid group action on $\mathcal{A}_3$}

Let us recall that we define the generators $Q_{13}$ as follows
\begin{align}\label{eq:Q13d}
Q_{13}=Q_{13d}&=\czR_{2}^{-1}Q_{12}\czR_{2}=\czR_{1}Q_{23}\czR_{1}^{-1}.
\end{align}
From the result of Proposition \ref{prop:CQ}, we know that $Q_{13}$ satisfies
\begin{align}
 Q_{13}
=\frac{Q_{1}Q_{3}+Q_{2}Q_{123}}{q+q^{-1}}-\frac{[Q_{12},Q_{23}]_q}{q^{2}-q^{-2}},
\end{align}
and is in the algebra $\mathcal{A}_3$ which is generated by $Q_1$, $Q_2$, $Q_3$, $Q_{12}$,
$Q_{23}$ and $Q_{123}$.  Now from \eqref{eq:Q13d} it is natural to consider the following
element which is analogous to
$Q_{13d}$:
\begin{align}\label{eq:Q13u}
Q_{13u}&=\czR_{2}Q_{12}\czR_{2}^{-1}=\czR_{1}^{-1}Q_{23}\czR_{1}.
\end{align}
It has been shown in \cite{CrampeGaboriaudetal2020} that this element is also in
$\mathcal{A}_3$ since it can be obtained as
\begin{align}
Q_{13u}
=\frac{Q_{1}Q_{3}+Q_{2}Q_{123}}{q+q^{-1}}-\frac{[Q_{23},Q_{12}]_q}{q^{2}-q^{-2}}.
\end{align}
The labels $u$ and $d$ added on the Casimir elements $Q_{13d}$ and $Q_{13u}$ stand for
\textit{up} and \textit{down}.  These names come from the form of their image in
$\Sk_{\,iq^{1/2}}(\Sigma_{0,4})$ given in Corollary \ref{co:CQ}:
\begin{subequations}
\begin{align}
&\begin{aligned}
\begin{tikzpicture}
 \punctiii \darc{2} \darc{0} \bendu{0.25}{1.75} \bendu{-0.25}{2.25} \draw (2.5,0)
 node[label=right:{$= \phi(Q_{13u}),$}]{};
\end{tikzpicture}
\end{aligned}\\
&\begin{aligned}
\begin{tikzpicture}
 \punctiii \uarc{2} \uarc{0} \bendd{0.25}{1.75} \bendd{-0.25}{2.25} \draw (2.5,0)
 node[label=right:{$= \phi(Q_{13d}).$}]{};
\end{tikzpicture}
\end{aligned}
\end{align}
\end{subequations}
This procedure of obtaining additional elements of $\mathcal{A}_3$ by conjugations of
braided $R$-matrices can be described by an automorphism action. Let us define the
following automorphisms of $\mathcal{A}_3$ denoted $\Psi_{s_i}$ and $\Psi_{s_i^{-1}}$ by
\begin{align}
 \Psi_{s_i}(X)= \czR_i X \czR_i^{-1} \qquad\text{and}
 \qquad \Psi_{s_i^{-1}}(X)= \czR_i^{-1} X \czR_i= \Psi_{s_i}^{-1}(X),
\end{align}
for $i=1,2$ and $X\in \mathcal{A}_3$, The previous maps are well-defined since the images
of the generators of $\mathcal{A}_3$ are precisely in $\mathcal{A}_3$ (and not in its
completion).  Indeed, by direct computations making use of the explicit form \eqref{eq:uR}
of the universal $R$-matrix and the commutation relations of $U_q(\sl_2)$, one gets
\begin{gather}\label{psi1}
\begin{gathered}
\Psi_{s_1}(Q_1)=Q_2, \qquad
\Psi_{s_1}(Q_2)=Q_1, \qquad
\Psi_{s_1}(Q_3)=Q_3, \qquad
\Psi_{s_1}(Q_{123})=Q_{123},\\
\Psi_{s_1}(Q_{12})=Q_{12}, \qquad
\Psi_{s_1}(Q_{23})=Q_{13d}.
\end{gathered}
\end{gather}
and
\begin{gather}\label{psi2}
\begin{gathered}
\Psi_{s_2}(Q_1)=Q_1, \qquad
\Psi_{s_2}(Q_2)=Q_3, \qquad
\Psi_{s_2}(Q_3)=Q_2, \qquad
\Psi_{s_2}(Q_{123})=Q_{123},\\
\Psi_{s_2}(Q_{12})=Q_{13u}, \qquad
\Psi_{s_2}(Q_{23})=Q_{23}.
\end{gathered}
\end{gather}
We obtain similarly the actions of $\Psi_{s_i^{-1}}$ on the generators of $\mathcal{A}_3$.

Since the braided $R$-matrix satisfies the braided Yang--Baxter equation \eqref{bYBE},
we can show that the defining relations \eqref{eq:braid_relations} of the braid group
$B_3$ are reproduced
\begin{subequations}\label{eq:psi}
\begin{align}
 \Psi_{s_1}\circ\Psi_{s_2}\circ\Psi_{s_1}&=\Psi_{s_2}\circ\Psi_{s_1}\circ\Psi_{s_2},\\
 \Psi_{s_i}\circ\Psi_{s_i^{-1}}&=\Psi_{s_i^{-1}}\circ \Psi_{s_i}=id.
\end{align}
\end{subequations}
We extend the automorphisms $\Psi_S$ to any $S\in B_3$ by
\begin{equation}\label{PsiS}
\Psi_S(X) = (\Psi_{g_1}\circ\Psi_{g_2}\circ\dots\circ\Psi_{g_\ell})(X),
\end{equation}
where $S$ is decomposed as $S=g_1g_2\dots g_\ell$ and
$g_i\in\{s_1,s_2,s_1^{-1},s_2^{-1}\}$.  Note that the map \eqref{PsiS} does not depend on
the choice of the decomposition of $S$ due to \eqref{eq:psi}.
\begin{rmk}
The realization of the braid group given by $\Psi_S$ is not faithful. For example, one can
verify that $\Psi_{(s_1s_2)^{3}}=id$. This is checked to be true on the intermediate
Casimir elements by making repeated use of \eqref{psi1}-\eqref{psi2}. It follows that it
is also true for any polynomial in those elements. Moreover, some elements of
$\mathcal{A}_3$ have additional stabilizers, \textit{e.g.}
\begin{subequations}
\begin{align}
 \Psi_{s_1s_1}(Q_1)&=\czR_1^{-1}\czR_1^{-1} Q_1 \czR_1\czR_1=\czR_1^{-1}Q_2\czR_1=Q_1,\\
 \Psi_{s_2}(Q_{23})&=Q_{23}.
\end{align}
\end{subequations}
Identifying stabilizers of the braid group action on elements of $\mathcal{A}_3$ is easy
to do but giving an exhaustive list is harder.
\end{rmk}
\begin{rmk}
It was shown in \cite{CrampeVinetetal2020b} how such a braid group action translates to
the $q\to-1$ limit. This limit of the Askey--Wilson algebra is referred to as the
Bannai--Ito algebra. In that case, the $B_3$ braid group action simplifies to an action of
the $S_3$ symmetric group. It is possible to study more generally the action of the
$S_n$ symmetric group on the higher rank Bannai-Ito algebra $B(n)$.
\end{rmk}

\subsection{Half Dehn twists and the braid group action on
$\Sk_{\,iq^{1/2}}(\Sigma_{0,4})$}

We now present a $B_3$ group action on the Kauffman bracket skein algebra
$\Sk_{\,iq^{1/2}}(\Sigma_{0,4})$, denoted
$\psi_{S}:\Sk_{\,iq^{1/2}}(\Sigma_{0,4})\to\Sk_{\,iq^{1/2}}(\Sigma_{0,4})$, with $S\in
B_3$.  The braid group action rotates the placement of the punctures with respect to each
other.

Here is how it goes. First, the actions $\psi_{s_i}$ and $\psi_{s_i^{-1}}$ on
$\Sk_{\,iq^{1/2}}(\Sigma_{0,4})$ are defined by the so-called half Dehn twists
\cite{FarbMargalit2012, Hikami2019}.  The four generators of $B_3$ act as
\begin{align}\label{eq:B3actions}
\begin{aligned}
\begin{tikzpicture}\punctiii \arrur{0} \arrdl{0}
\draw (-0.5,0) node[label=left:{$\psi_{s_1^{\phantom{-1}}}=$}]{};
\end{tikzpicture},\\
\begin{tikzpicture}\punctiii \arrul{0} \arrdr{0}
\draw (-0.5,0) node[label=left:{$\psi_{s_1^{-1}}=$}]{};
\end{tikzpicture},
\end{aligned}\hspace{6em}
\begin{aligned}
\begin{tikzpicture}\punctiii \arrur{1} \arrdl{1}
\draw (-0.5,0) node[label=left:{$\psi_{s_2^{\phantom{-1}}}=$}]{};
\end{tikzpicture},\\
\begin{tikzpicture}\punctiii \arrul{1} \arrdr{1}
\draw (-0.5,0) node[label=left:{$\psi_{s_2^{-1}}=$}]{};
\end{tikzpicture},
\end{aligned}
\end{align}
where any framed link gets deformed continuously without crossing the punctures as the
rotations happen.  For example, one gets
\begin{align}
\begin{aligned}
\psi_{s_2^{-1}}(\mathbb{A}_{12})=\psi_{s_2^{-1}}\left(
\begin{aligned}
\begin{tikzpicture}
\punctiii \larc{0} \updown{0}{1} \rarc{1}
\end{tikzpicture}
\end{aligned}
\right)&=\left(
\begin{aligned}
\begin{tikzpicture}
\punctiii \larc{0} \updown{0}{1} \rarc{1} \arrul{1} \arrdr{1}
\end{tikzpicture}
\end{aligned}
\right)\\
&=\left(
\begin{aligned}
\begin{tikzpicture}
\punctiii \uarc{2} \uarc{0} \bendd{0.25}{1.75} \bendd{-0.25}{2.25}
\end{tikzpicture}
\end{aligned}\right)=\mathbb{A}_{13},
\end{aligned}
\end{align}
and
\begin{align}
\begin{aligned}
\psi_{s_2^{\phantom{-1}}}(\mathbb{A}_{23})=
\psi_{s_2^{\phantom{-1}}}\left(
\begin{aligned}
\begin{tikzpicture} \punctiii \larc{1}\updown{1}{2}\rarc{2} \end{tikzpicture}
\end{aligned}
\right)
&=\left(
\begin{aligned}
\begin{tikzpicture}
\punctiii \larc{1}\updown{1}{2}\rarc{2}\arrur{1}\arrdl{1}
\end{tikzpicture}
\end{aligned}
\right)\\&=\left(
\begin{aligned}
\begin{tikzpicture} \punctiii \larc{1}\updown{1}{2}\rarc{2} \end{tikzpicture}
\end{aligned}
\right)
=\mathbb{A}_{23}.
\end{aligned}
\end{align}
\begin{prop}
The actions $\psi_g$ for $g\in \{s_1,s_2,s_1^{-1},s_2^{-1}\}$ are automorphisms of
$\Sk_{\,iq^{1/2}}(\Sigma_{0,4})$.
\end{prop}
\proof For any $X$, $Y\in\Sk_{\,iq^{1/2}}(\Sigma_{0,4})$ and $g\in
\{s_1,s_2,s_1^{-1},s_2^{-1}\}$, one understands that
\begin{align}\label{eq:homProperty}
 \psi_g(X\cdot Y)=\psi_g(X)\cdot\psi_g(Y).
\end{align}
Indeed, from the way they were defined, the rotations do not add or change crossings.
Thus, the Kauffman bracket relations \eqref{eq:Kauf} that one makes use of to ``simplify
the crossings'' of a given product are unchanged under these rotations. Since the
rotations are also defined in order to avoid links crossing punctures, the topological
properties (such as which punctures are circled by which links) are preserved. Hence the
action $\psi_g$ is a homomorphism.  Moreover $\psi_g$ is an endomorphism because links in
$\Sk_{\,iq^{1/2}}(\Sigma_{0,4})$ are mapped to other links in
$\Sk_{\,iq^{1/2}}(\Sigma_{0,4})$, and it is invertible, as rotations can be inverted, thus
$\psi_g$ is an automorphism.  \endproof

Let $S=g_1g_2\dots g_\ell \in B_3$ be a decomposition of an element of the braid group on
three strands with $g_i\in\{s_1,s_2,s_1^{-1},s_2^{-1}\}$.  We define the automorphism
$\psi_S$ as follows:
\begin{equation}\label{psiS}
 \psi_S(X) = (\psi_{g_1}\circ\psi_{g_2}\circ\dots\circ\psi_{g_\ell})(X).
\end{equation}
We use also the definition $\psi_1=id$.
The previous map \eqref{psiS} does not depend on the choice of the decomposition of $S$.
Indeed, it is straightforward to check that the defining relations of the braid group
\eqref{eq:braid_relations} are verified on the generators. By the homomorphism property
\eqref{eq:homProperty}, it follows that these braid relations are verified for any element
of the Kauffman bracket skein module $\Sk_{\,iq^{1/2}}(\Sigma_{0,4})$.
\begin{rmk}
More visually complicated loops can always be created by further ``twisting'' the loops.
For example,
\begin{align}\label{eq:circlingloop}
\psi_{(s_1^{-1})^{4}}(\mathbb{A}_{23})=
\begin{aligned}\begin{tikzpicture}
\punct{0}\punct{1}\punct{2.25}
\draw[fill=black, fill opacity=0.075] (0.75cm,0) arc (180:0:0.25cm) to [bend left=90]
 (-0.5cm,0) to [bend left=90] (1.75cm,0) to [bend left=90] (-1.cm,0)
 to [bend left=90] (2.5cm,0) arc (0:-180:0.25cm) to [bend right=90] (-0.75cm,0) to
 [bend right=90] (1.5cm,0) to [bend right=90] (-0.25cm,0) to [bend right=90] (0.75cm,0);
\end{tikzpicture}\end{aligned}
\end{align}
is a more complicated analog of $\mathbb{A}_{23}$. The shadow filling the inside of the
loop is there to guide the eyes of the reader.  These have also been studied in
\cite{BakshiMukherjeeetal2018}.
\end{rmk}
\begin{rmk}
Let us remark that in \cite{Hikami2019}, the author considers a similar braid group action
by half Dehn twists on the Kauffman bracket skein algebra of the four-punctured sphere. In
that paper, it is shown that the group $SL(2;\mathbb{Z})$ acts on the DAHA of type
$(C_1^{\vee}C_1)$ through conjugations. Furthermore, the Artin braid group $B_3$ action on
the Kauffman bracket skein algebra can be seen as a translation of this $SL(2;\mathbb{Z})$
action. We also note that Terwilliger had presented a $B_3$ action on both the universal
Askey--Wilson algebra and the universal DAHA of type $(C_1^{\vee},C_1)$
\cite{Terwilliger2013}.
\end{rmk}

\subsection{Connection between both braid actions}

The following proposition establishes the connections between both braid group actions
presented above.
\begin{prop}
The following diagram of isomorphisms
\begin{center}
\begin{tikzpicture}[scale=0.5]
\node (skg) at (0,5) {$\mathcal{A}_3$};
\node (skd) at (5,5) {$\mathcal{A}_3$};
\node (Ag) at (0,0) {$\Sk$};
\node (Ad) at (5,0) {$\Sk$};
\node at (2.5,5.5) {$\Psi_S$};
\node at (2.5,0.5) {$\psi_S$};
\node at (0.5,2.5) {$\phi$};
\node at (5.5,2.5) {$\phi$};
\draw [thick,->] (skg) -- (skd);
\draw [thick,->] (skg) -- (Ag);
\draw [thick,->] (Ag) -- (Ad);
\draw [thick,->] (skd) -- (Ad);
\end{tikzpicture}
\end{center}
is commutative for any $S\in B_3$. Here we used the shortened notation
$\Sk\equiv\Sk_{\,iq^{1/2}}(\Sigma_{0,4})$. The isomorphisms $\phi$, $\Psi_S$ and $\psi_S$
are given in \eqref{eq:cSQ}, \eqref{PsiS} and \eqref{psiS}, respectively.
\end{prop}
\proof
We can show that this diagram is commutative for all the generators of $\mathcal A$ and
for any $S=s_i$ or $S=s_i^{-1}$. For example:
\begin{align}
 \phi\circ \Psi_{s_1}( Q_1)
 = \phi(Q_2)
 = \mathbb{A}_2
 = \psi_{s_1} ( \mathbb{A}_1)
 = \psi_{s_1} \circ \phi (Q_1).
\end{align}
A more complicated example is
\begin{align}
 \phi\circ \Psi_{s_2}( Q_{12} )
 = \phi(Q_{13u})
 = \begin{aligned}\begin{tikzpicture}
    \punctiii \darc{2} \darc{0} \bendu{0.25}{1.75} \bendu{-0.25}{2.25}
   \end{tikzpicture}\end{aligned}
 = \psi_{s_2}(\mathbb{A}_{12})
 = \psi_{s_2}\circ \phi ( Q_{12} ).
\end{align}
Since all the maps of the diagram are homomorphisms, the commutativity of the diagram on
the generators of $\mathcal A$ is enough to prove the proposition for any $S\in B_3$.
\endproof

The commutativity of this diagram allows us to identify the conjugation by the braided
$R$-matrix for $\mathcal{A}_3$ as half Dehn twists around the punctures of
$\Sk_{\,iq^{1/2}}(\Sigma_{0,4})$.  In addition, we can identify easily the elements of the
algebra $\mathcal{A}_3$ obtained as an image by $\Psi_S$ with a link of
$\Sk_{\,iq^{1/2}}(\Sigma_{0,4})$.

\section{Towards a higher rank $\saw(n)$ algebra}\label{sec:higher}

Some natural generalizations of the different algebras have previously been introduced and
studied:
\begin{itemize}
\item the generalized Askey--Wilson algebra $\mathbf{aw}(n)$ is the
 algebra generated by $\{ C_I\ |\  I \subset \{1,2,\dots,n\}\}$ subject to the relations
 introduced in Theorems 3.1 and 3.2 of \cite{DeClercq2019};
\item the algebra $\mathcal{A}_n$ is the subalgebra of $U_q(\sl_2)^{\otimes n}$ generated
 by all the intermediate Casimir elements $\{ Q_I\ |\  I \subset \{1,2,\dots,n\}\}$
 obtained by the repeated action of the coproduct of $U_q(\sl_2)$;
\item the centralizer $\mC_n$ is defined by
 \begin{equation}\label{eq:defCn}
 \mC_n=\{ X \in U_q(\sl_2)^{\otimes n}\ \big| \ [ \Delta^{(n-1)}(x), X]=0\ ,
 \ \  \forall x\in U_q(\sl_2) \}
 \end{equation}
 where $\Delta^{(n)}=(\Delta^{(n-1)}\otimes id)\Delta$ and $\Delta^{(1)}=\Delta$;
\item the algebra $Sk_{\theta}(\Sigma_{0,n+1})$ is the Kauffman bracket skein algebra
 associated to the $(n+1)$-punctured sphere $\Sigma_{0,n+1}$ \cite{CookeLacabanne2020}.
 Let us now associate to each set $I\subseteq[1;n]\equiv\{1,2,\dots,n\}$ a `simple' loop
 $\mathbb{A}_I$ of $Sk_{\theta}(\Sigma_{0,n+1})$.  We write a set $I$ as $I=I_1\cup
 I_2\cup\dots\cup I_\ell$, where $I_i$ are sets of consecutive integers and then we define
 the `simple' loop $\mathbb{A}_I$ as:
 \begin{align}
 \mathbb{A}_I&=\qquad
 \begin{aligned}
 \begin{tikzpicture}
 \draw (0,0) node[]{$\dots$};
 \draw (1,0) node[label=above:{\scriptsize$I_1$}]{$\times$};
 \draw (1.5,0) node[]{$...$};
 \draw (2,0) node[label=above:{\scriptsize$I_2$}]{$\times$};
 \draw (3,0) node[]{$\dots$};
 \draw (4,0) node[label=above:{\scriptsize$I_\ell$}]{$\times$};
  \draw (5,0) node[]{$\dots$};
 \uarc{1} \darc{1.5} \uarc{2} \uarc{4} \bendd{0.75}{4.25}
 \draw (2.25cm,0) arc (180:270:0.25cm);
 \draw (3.75cm,0) arc (0:-90:0.25cm);
 \end{tikzpicture}
 \end{aligned}
 \end{align}
 These simple loops do not bend around, unlike \eqref{eq:circlingloop}. They are only
 extending in the lower half of the plane. In particular, for $I=\{i,i+1,\dots ,j\}$, a
 set of consecutive integers, one gets
 \begin{align}
 \mathbb{A}_I=\qquad\left(
 \begin{aligned}
 \begin{tikzpicture}
 \draw (0,0) node[label=below:{\scriptsize$1$}]{$\times$}; \draw (1,0) node[]{$\dots$};
 \draw (2,0) node[label=below:{\scriptsize$i$}]{$\times$}; \draw (2.5,0) node[]{$\dots$};
 \draw (3,0) node[label=below:{\scriptsize$j$}]{$\times$};
 \draw (4,0) node[]{$\dots$}; \draw (5,0) node[label=below:{\scriptsize$n$}]{$\times$};
 \larc{2} \rarc{3} \updown{2}{3}
 \end{tikzpicture}
 \end{aligned}\right)\qquad=\qquad\left(
 \begin{aligned}
 \begin{tikzpicture}
 \draw (0,0) node[]{$\dots$}; \draw (1,0) node[label=below:{\scriptsize$I$}]{$\times$};
 \draw (2,0) node[]{$\dots$};
 \larc{0.75} \rarc{1.25} \updown{0.75}{1.25};
 \end{tikzpicture}
 \end{aligned}\right)
 \end{align}
 \end{itemize}

What is lacking in the previous list is the generalization $\saw(n)$ of the algebra
$\saw(3)$. Such a generalization would provide a description of the algebra
$\mathcal{A}_n$ in terms of generators and relation.  We know that $\saw(n)$ will be a
quotient of the algebra $\mathbf{aw}(n)$ by relation(s) of the type \eqref{eq:awt3cas1},
with some Casimir elements to be determined.
\textit{We conjecture that the map $\phi_n$ from $\saw(n)$ to
$Sk_{\theta}(\Sigma_{0,n+1})$ which sends $Q_I$ to $\mathbb{A}_I$ is an isomorphism}
\footnote{During the preparation of this paper,
the authors have been informed by J. Cooke that a similar idea was pursued in an upcoming
publication \cite{CookeLacabanne2020}.}.

Let us mention that there also exist generalizations in the non-deformed case ($q=1$ and
$q=-1$) of the Askey--Wilson algebra: these are respectively called the ``higher rank
Racah algebra'' introduced in \cite{DeBieGenestetal2017} as well as the ``higher rank
Bannai--Ito algebra'' introduced in \cite{DeBieGenestetal2016}.

In the remainder, we give different indications regarding ways to define $\saw(n)$.

\subsection{Punctures on a sphere and a coassociative homomorphism of Kauffman bracket
skein modules}

Recall we had highlighted that the punctures of the sphere were related to the tensor
product factors. Additionally, a loop encircling a puncture is associated to some
intermediate Casimir element with non-trivial factors in the tensor product factor
corresponding to the puncture.

Further recall that the coproduct $\Delta$ acts as an algebra morphism from $U_q(\sl_2)$
to $U_q(\sl_2)^{\otimes2}$.  One can define an action of the coproduct on any
$i$\textsuperscript{th} factor of a tensor product: for any $X\in U_q(\sl_2)^{\otimes n}$,
we define $\Delta_{i}:U_q(\sl_2)^{\otimes n}\to U_q(\sl_2)^{\otimes(n+1)}$ as:
\begin{align}
 \Delta_{i}(x)=\left(1^{\otimes(i-1)}\otimes\Delta\otimes1^{\otimes(n-i)}\right)(X).
\end{align}
Now in $U_q(\sl_2)^{\otimes3}$ some intermediate Casimir elements are related to each
other by the coproduct, such as $Q_1$ and $Q_{12}$:
\begin{align}\label{eq:L1L12Delta}
 \Delta_1(Q_1)=(\Delta\otimes1\otimes1)Q_1=Q_{12}\otimes1.
\end{align}
This relation between $Q_1$ and $Q_{12}$ appears in the framed links picture as well.

More precisely, $\Delta_i$ has an analog, the $\delta_i$ morphism, which acts on a
single puncture $i$ by creating another puncture next to it.  If the puncture $i$ is
enclosed in a loop, the created puncture is also enclosed in the same loop. The example
\eqref{eq:L1L12Delta} is illustrated as follows:
\begin{align}\label{eq:C1C12delta}
\begin{aligned}
\delta_{1}\mathbb{A}_1
=\delta_1\left(\begin{aligned}\begin{tikzpicture}
\punctiii \larc{0} \rarc{0}
\end{tikzpicture}\end{aligned}\right)
=\delta_1&\left(\begin{aligned}\begin{tikzpicture}
\draw (0.5,0) node{$\times$};
\draw (2,0) node{$\times$};
\draw (3,0) node{$\times$};
\larc{0} \rarc{1} \updown{0}{1}
\end{tikzpicture}\end{aligned}\right)\\
=&\left(\begin{aligned}\begin{tikzpicture}
\punctiv \larc{0} \rarc{1} \updown{0}{1}
\end{tikzpicture}\end{aligned}\right)=\mathbb{A}_{12}\in\Sk_{\theta}(\Sigma_{0,5})
\end{aligned}
\end{align}
This $\delta_i$ is a Kauffman bracket skein module coassociative algebra homomorphism.
It provides embeddings of $\Sk_{\theta}(\Sigma_{0,n})\to\Sk_{\theta}(\Sigma_{0,n+1})$.
This can be seen as the commutativity of the following diagram:
\begin{center}
\begin{tikzpicture}[scale=0.5]
\node (skg) at (0,5) {$\mathcal{A}_n$};
\node (skd) at (8,5) {$\mathcal{A}_{n+1}$};
\node (Ag) at (0,0) {$\Sk_{\theta}(\Sigma_{0,n+1})$};
\node (Ad) at (8,0) {$\Sk_{\theta}(\Sigma_{0,n+2})$};
\node at (3.5,5.5) {$\Delta_i$};
\node at (3.5,0.5) {$\delta_i$};
\node at (0.6,2.5) {$\phi_n$};
\node at (8.95,2.5) {$\phi_{n+1}$};
\draw [thick,->] (skg) -- (skd);
\draw [thick,->] (skg) -- (Ag);
\draw [thick,->] (Ag) -- (Ad);
\draw [thick,->] (skd) -- (Ad);
\end{tikzpicture}
\end{center}

\subsection{A crossing index}

The defining algebra relations of $\Sk_{\,iq^{1/2}}(\Sigma_{0,4})$
\eqref{eq:aw3}--\eqref{eq:awt3cas1} (see Proposition \ref{pr:as}) can be classified in
three types.  The relations always involve two generators, whose product, commutator or
$q$-commutator is reexpressed in terms of other generators. Now imagine we draw both
generators simultaneously in a framed links diagram (as if we were to multiply them).
Some crossings will appear if the two generators don't commute.
\begin{defi}
The \textit{crossing index} is defined as the minimal number of crossings that appear in a
framed link diagram no matter how the generators are drawn.
\end{defi}
The relations \eqref{eq:aw3}--\eqref{eq:awt3cas1} can be classified in terms of the
crossing index as follows:
\begin{itemize}
\item If the generators can be drawn simultaneously in such a way that the loops have no
crossings (crossing index of $0$), they will commute (for example, this is the case for
any central element $Q_1,Q_2,Q_3,Q_{123}$ multiplied with any other generator).
\item If the generators can be drawn in such a way that their minimum number of crossings
is two (crossing index of $2$), linear $q$-commutation relations of $\mathbf{aw}(3)$-type
will be obtained, such as relations \eqref{eq:aw3}.
\item If the generators have a crossing index of $4$, such as
\begin{align}\label{}
\phi(Q_{13u}Q_{13d})=
\begin{aligned}\begin{tikzpicture}
 \punctiii 
 \draw (2.5,0) node{,};
 \draw (-0.25cm,0.2cm) to [bend left=90] (2.25cm,0.2cm);
 \draw (0.25cm,0.2cm) to [bend left=90] (1.75cm,0.2cm);
 \draw (-0.2cm,-0.05cm) arc (190:350:0.2cm);
 \draw (1.8cm,-0.05cm) arc (190:350:0.2cm);
 \uarc{2} \uarc{0} \bendd{0.25}{1.75} \bendd{-0.25}{2.25}
\end{tikzpicture}\end{aligned}
\end{align}
higher order relations of the type
\eqref{eq:awt3cas1} will be obtained.
\end{itemize}
This crossing index proves useful for the analysis of the higher rank generalizations of
$\saw(3)$.

\subsection{The algebras $\mathbf{aw}(n)$ and $Sk_{\theta}(\Sigma_{0,n+1})$}

As mentioned previously, the algebra $\mathbf{aw}(n)$ is generated by $C_I$ with
$I\subseteq[1;n]$ and subject to the relations of Proposition 3.1 of \cite{DeClercq2019}.
We can show by using the action of the morphism $\delta_i$
that we have an homomorphism from $\mathbf{aw}(n)$ to $Sk_{\theta}(\Sigma_{0,n+1})$.
Moreover, we can show that all the relations of Proposition 3.1 of \cite{DeClercq2019}
correspond to the product of two simple loops with crossing index $2$.  We believe that
the relations in \cite{DeClercq2019} exhaust all possibilities of relations involving the
product of simple loops with crossing index 2.  We conjecture also that the above
mentioned homomorphism is surjective (but it is certainly not injective, even for the case
$n=3$).  The description of the kernel would involve products of links with a crossing
index strictly greater than 2.  The complete description of this kernel would lead to the
definition of $\saw(n)$ and give an algebraic description of $\mathcal A_n$ and
$Sk_{\theta}(\Sigma_{0,n+1})$.

The study of $\saw(n)$ should be guided by the intuition gained from the
framed links picture.  To illustrate the type of insight we can gain, let us
efficiently summarize some of the results of \cite{PostWalter2017}. In this paper, the
authors study the intermediate Casimir elements in $U_q(\sl_2)^{\otimes 4}$ and
introduce an involution $I$ of the algebra as well as ``involuted'' generators $IQ_{13}$
and $IQ_{24}$ satisfying
\begin{equation}
 [Q_{13},IQ_{24}] =0 , \qquad\text{and}\qquad [IQ_{13},Q_{24}]=0.
\end{equation}
That these generators commute becomes evident when we rewrite (following our definitions)
$IQ_{24}=Q_{24u}$, $IQ_{13}=Q_{13u}$, and then draw the corresponding links. Indeed, the
products
\begin{subequations}
\begin{align}
\phi(Q_{13d}Q_{24u})&=\begin{aligned}
\begin{tikzpicture}
\punctiv \uarc{0} \darc{1} \uarc{2} \darc{3} \bendd{-0.25}{2.25} \bendd{0.25}{1.75}
\bendu{0.75}{3.25} \bendu{1.25}{2.75}
\end{tikzpicture}
\end{aligned}=\phi(Q_{24u}Q_{13d}),\\[-1em]
\phi(Q_{13u}Q_{24d})&=\begin{aligned}
\begin{tikzpicture}
\punctiv \darc{0} \uarc{1} \darc{2} \uarc{3} \bendu{-0.25}{2.25} \bendu{0.25}{1.75}
\bendd{0.75}{3.25} \bendd{1.25}{2.75}
\end{tikzpicture}
\end{aligned}=\phi(Q_{24d}Q_{13u}),
\end{align}
\end{subequations}
have $0$ crossing hence $[Q_{13d},Q_{24u}]=0$ and $[Q_{13u},Q_{24d}]=0$.

What about the product of terms like $Q_{13d}$ and $Q_{24d}$ ?  This calculation has never
appeared in the papers mentioned above because it has a crossing number of $4$:
\begin{align}
\phi(Q_{13d}Q_{24d})&=
\begin{aligned}
\begin{tikzpicture}
\punctiv \uarc{0} \uarc{1} \uarc{2} \uarc{3}
\draw (0.25cm,0) arc (180:270:0.25cm);
\draw (0.5cm,-0.25cm) -- (0.7cm,-0.25cm);
\draw (0.8cm,-0.25cm) -- (1.2cm,-0.25cm);
\draw (1.3cm,-0.25cm) -- (1.5cm,-0.25cm);
\draw (1.5cm,-0.25cm) arc (270:360:0.25cm);
\draw (-0.25cm,0cm) -- (-0.25cm,-0.25cm);
\draw (-0.25cm,-0.25cm) arc (180:270:0.25cm);
\draw (0cm,-0.5cm) -- (0.7cm,-0.5cm);
\draw (0.8cm,-0.5cm) -- (1.2cm,-0.5cm);
\draw (1.3cm,-0.5cm) -- (2cm,-0.5cm);
\draw (2cm,-0.5cm) arc (270:360:0.25cm);
\draw (2.25cm,-0.25cm) -- (2.25cm,0.cm);
\draw (1.25cm,0cm) -- (1.25cm,-0.5cm);
\draw (1.25cm,-0.5cm) arc (180:270:0.25cm);
\draw (1.5cm,-0.75cm) -- (2.5cm,-0.75cm);
\draw (2.5cm,-0.75cm) arc (270:360:0.25cm);
\draw (2.75cm,-0.5cm) -- (2.75cm,0cm);
\draw (0.75cm,0cm) -- (0.75cm,-0.75cm);
\draw (0.75cm,-0.75cm) arc (180:270:0.25cm);
\draw (1.cm,-1cm) -- (3cm,-1cm);
\draw (3cm,-1cm) arc (270:360:0.25cm);
\draw (3.25cm,-0.75cm) -- (3.25cm,0cm);
\end{tikzpicture}
\end{aligned}
\end{align}
Remarkably, this calculation can be effected in $\Sk_{\,iq^{1/2}}(\Sigma_{0,5})$ using the
conjectured morphism. One writes the $Q_I$ in terms of $\mathbb{A}_I$, computes using the
skein relations of $\Sk_{\,iq^{1/2}}(\Sigma_{0,5})$, then reexpresses all $\mathbb{A}_I$
in terms of $Q_I$. This yields the following results
\begin{align}
&Q_{13d} \, Q_{24d} = q^2 Q_{14d}Q_{23} + q^{-2} Q_{12}Q_{34} - q(Q_{14d}Q_2Q_3 +
Q_{23}Q_1Q_4) - q^{-1}(Q_{12}Q_3Q_4 + Q_{34}Q_1Q_2) \nonumber \\
&\hspace{3em}
- (q+q^{-1})Q_{1234} + Q_1Q_2Q_3Q_4 + Q_1Q_{234} + Q_2Q_{134d} + Q_3Q_{124d}+ Q_4Q_{123}
\intertext{and}
&Q_{24d} \, Q_{13d} = q^{-2} Q_{14d}Q_{23} + q^{2} Q_{12}Q_{34} - q^{-1}(Q_{14d}Q_2Q_3 +
Q_{23}Q_1Q_4) - q(Q_{12}Q_3Q_4 + Q_{34}Q_1Q_2) \nonumber \\
&\hspace{3em}
 - (q+q^{-1})Q_{1234} + Q_1Q_2Q_3Q_4 + Q_1Q_{234} + Q_2Q_{134d} + Q_3Q_{124d}+ Q_4Q_{123}.
\end{align}
These have been checked to hold in $U_q(\sl_2)^{\otimes4}$.

Let us also mention that the action of the braid group $B_3$ can be generalized to the
action of $B_n$ on $Sk_{\theta}(\Sigma_{0,n+1})$ and $\mathcal{A}_n$. This might turn out
useful for proving results in the future.

\section{Conclusion}\label{sec:concl}

Three objectives were principally pursued in this paper. The first aimed to review the
different avatars of the Askey--Wilson algebra and to clarify the relations between them.
Among those algebras, we focused on two and presented novel results related to these
cases; this was the second main goal.  The Special Zhedanov algebra
$\szh(m_1,m_2,m_3,m_4)$ was obtained from (a quotient of) the reflection algebra by
setting the Sklyanin determinant to an appropriate value; its $W(D_4)$ symmetry was
exhibited in addition.  The Special Askey--Wilson algebra $\saw(3)$, a universal analogue
of $\szh(m_1,m_2,m_3,m_4)$, was shown to be isomorphic to the algebra $\mathcal{A}_3$ that
emerges from the Racah problem of $U_q(\sl_2)$ and also to the Kauffman bracket skein
algebra of the four-punctured sphere $\Sk_{\,iq^{1/2}}(\Sigma_{0,4})$.  An injective
homomorphism between $\mathcal{A}_3$ and the centralizer $\mC_3$ of $U_q(\sl_2)$ in its
threefold tensor product was stated and proved.  Actions of the braid group on both
$\Sk_{\,iq^{1/2}}(\Sigma_{0,4})$ (through half Dehn twists) and $\mathcal{A}_3$ (through
conjugation by braided $R$-matrices) were illustrated and shown to be compatible.  The
third main objective was to discuss the generalization of $\saw(3)$ to $\saw(n)$. To that
end, we emphasized the diagrammatic approach, defined a crossing index, and revisited the
results of \cite{PostWalter2017, DeClercq2019} in a unified fashion.

Let us conclude with more remarks regarding generalizations of Askey--Wilson algebras. It
would certainly be desirable to return to Zhedanov's original quest and to determine
directly from the multivariate Askey--Wilson polynomials (of Tratnik type)
\cite{GasperRahman2005} the algebra that encapsulates their bispectral properties. Steps
have been carried out \cite{Iliev2011, BaseilhacMartin2018, DeBieDeClercq2020} but this
should be completed. A definite higher rank generalization of the Zhedanov algebra will
emerge, whose quotients and central extensions could then be examined and should connect
to various fields in mathematics and physics.  Considering higher rank Lie algebras
$\mathfrak{g}$ instead of $\sl_2$ is another avenue that should be explored. The
centralizer of the diagonal action of $U_q(\mathfrak{g})$ in the $n$-fold tensor product
$U_q(\mathfrak{g})^{\otimes n}$, or the algebra generated by all the intermediate Casimir
elements of $\mathfrak{g}$ in the associated Racah problem should be studied.  Connections
with a generalization of $Sk_{\theta}(\Sigma_{0,n})$ to punctured manifolds of higher
genera would be worth investigating (see also \cite{CookeSamuelson2020}). We may also
wonder whether the braided universal $R$-matrix of $U_q(\mathfrak{g})$ plays a role in
this context.  Furthermore, the truncated reflection algebra presented in Section
\ref{sec:reflection} provides a natural framework to obtain generalizations of Zhedanov
algebras.  Different possibilities are here conceivable. One could consider more general
truncations of the reflection algebra.  This type of generalization has been already
studied in \cite{BaseilhacKoizumi2005a} and has been associated to quotients of
$q$-Onsager algebras \footnote{The classical limit $q\to1$ leads to subalgebras of the
loop algebra of $\mathfrak{sl}_2$ and to quotients of the Onsager algebra by Davis
relations \cite{BaseilhacCrampe2019}.}. Connections with centralizers and/or skein
algebras remain to be examined. Another possibility with respect to truncated reflection
algebras is the following.  Instead of using the $R$-matrix associated to quantum affine
algebras, one could consider the $R$-matrix corresponding to Yangians.  In this case, a
particular truncation of the reflection algebra leads to the Hahn algebra, which is a
specialization of the Zhedanov algebra, see \cite{CrampeRagoucyetal2019}.  Other
truncations should provide interesting generalizations of this algebra.  Finally, the FRT
presentation of the reflection algebra associated to higher rank Lie algebras and
superalgebras is well-known.  For instance, the twisted Yangians
$Y^{\text{tw}}(\mathfrak{o}_n)$ and $Y^{\text{tw}}(\mathfrak{sp}_n)$
\cite{MolevNazarovetal1996} and the reflection algebra $\mathcal{B}(n,\ell)$
\cite{MolevRagoucy2002} correspond to subalgebras of the Yangian of $\mathfrak{sl}_n$.
Some $q$-deformations of these structures have been also studied previously
\cite{MolevRagoucyetal2003} and are related to the quantum affine algebra of
$\mathfrak{sl}_n$. Their truncations have yet to be scrutinized and should possess
interesting features \footnote{Such an approach has been pursued in the classical limit
$q\to1$ \cite{BaseilhacCrampeetal2019} to obtain generalizations of the so-called
classical Askey--Wilson algebra and are seen as subalgebras of the $\mathfrak{sl}_n$
Onsager algebra \cite{UglovIvanov1996}.}.  These ideas that we plan on pursuing in the
near future are indications that there is much lying ahead with respect to algebras of the
Askey--Wilson type and what they will reveal and lead to.


\subsection*{Acknowledgments}

Many thanks to Geoffroy Bergeron for long-drawn discussions.
We have also benefitted from exchanging with Pascal Baseilhac, Juliet Cooke, Hendrik De
Bie, Hadewijch De Clercq, Sarah Post, Paul Terwilliger and Alexei Zhedanov.
N. Crampé and L. Poulain d'Andecy are partially supported by Agence Nationale de la
Recherche Projet AHA ANR-18-CE40-0001.
L. Frappat is grateful to the Centre de Recherches Mathématiques (CRM) for hospitality and
support during his visit to Montreal in the course of this investigation.
J. Gaboriaud holds an Alexander-Graham-Bell scholarship from the Natural Sciences and
Engineering Research Council of Canada (NSERC).
The research of L. Vinet is funded in part by a Discovery Grant from NSERC.
\\[2em]

\appendix

\section{Classical limit and injectivity}\label{app:inj}

We provide an explicit description of the classical limit of the realization of $\saw(3)$
in $U_q(\sl_2)^{\otimes 3}$ in terms of polarized traces, and use it to prove the
injectivity of the map from $\saw(3)$ to the centralizer $\mC_3$. In this appendix, we
will work with the formal series version of $U_q(\sl_2)$ and reduce the proof of the
injectivity statement to one in the universal enveloping algebra $U(sl_2)$, where we can
use known results of classical invariant theory involving polarized traces.

\subsection{Polarised traces in $U(\sl_2)^{\otimes 3}$}

The algebra $U(\sl_2)$ is generated by elements $e_{ij}$, $i,j\in\{1,2\}$, with the
defining relations $[e_{ij},e_{kl}]=\delta_{jk}e_{il}-\delta_{li}e_{kj}$ and
$e_{11}+e_{22}=0$. To join up with the notations used in the paper for
$U_q(\sl_2)$, we set $E=e_{12}$, $F=e_{21}$ and
$H=\frac{1}{2}(e_{11}+e_{22})=e_{11}=-e_{22}$, and the relations become:
\begin{equation}\label{relsl2}
 [H,E]=E,\qquad [H,F]=-F,\qquad [E,F]=2H.
\end{equation}
In a tensor product $U(\sl_2)^{\otimes N}$, we denote the generators by $e_{ij}^{(a)}$,
where $a\in\{1,\dots,N\}$ indicates the corresponding factor in the tensor product. The
polarized traces are the following elements:
\begin{align}
 T^{(a_1,\dots,a_d)}=e_{i_2i_1}^{(a_1)}e_{i_3i_2}^{(a_2)}\dots e_{i_1i_d}^{(a_d)},\qquad
 a_1,\dots,a_d\in\{1,\dots,N\},
\end{align}
where the summation over repeated indices is understood. The specific combinations of
polarized traces that will appear are:
\begin{gather}\label{gen-Racah}
 k_1:=T^{(1,1)},\quad k_2:=T^{(2,2)},\quad
 k_3:=T^{(3,3)},\quad k_4:=k_1+k_2+k_3+2(T^{(1,2)}+T^{(2,3)}+T^{(1,3)}),\nonumber\\[0.5em]
 X:=k_1+k_2+2T^{(1,2)},\qquad Y:=k_2+k_3+2T^{(2,3)},\qquad Z:=[X,Y].
\end{gather}
By direct computation, we can show that the expression of $Z$ in terms of polarized traces
is $Z=-8T^{(1,2,3)}$.

\subsection{The algebra $U_{\alpha}(\sl_2)$}

In this appendix, we will work with the formal series version of the quantum group
$U_q(\sl_2)$. We consider a formal parameter $\alpha$. The algebra $U_{\alpha}(\sl_2)$ is,
as a vector space, the space $U(\sl_2)[[\alpha]]$ of all formal power series in $\alpha$
with coefficients in $U(\sl_2)$, and the multiplication is determined by the defining
relations of $U_q(\sl_2)$, see section \ref{sect:Uqsl2}, where $q$ is replaced by
$e^{\alpha}$ and $q^H$ is replaced by $e^{\alpha H}$. This results in the following
relations deforming (\ref{relsl2}):
\begin{equation}\label{relsl2def}
 [H,E]=E,\qquad [H,F]=-F,\qquad
 [E,F]=\frac{e^{2\alpha H}-e^{-2\alpha H}}{e^{\alpha}-e^{-\alpha}}.
\end{equation}
Similarly, the algebra $U_{\alpha}(\sl_2)^{\hat{\otimes} N}$ is the vector space
$U(\sl_2)^{\otimes N}[[\alpha]]$ of formal series with coefficients in
$U(\sl_2)^{\otimes N}$ and multiplication induced by the above relations in each factor.
The comultiplication of $U_{\alpha}(\sl_2)$ is naturally obtained from the
comultiplication given for $U_q(\sl_2)$.

Note that the limit $\alpha\to 0$ from $U_{\alpha}(\sl_2)$ yields the algebra
$U(\sl_2)$ and the comultiplication becomes the diagonal embedding.

\subsection{Reduction to $U(\sl_2)$}

We want to prove that the following elements
\begin{equation}\label{setQ}
 Q_{12}^{i} Q_{23}^{j} Q_{13}^{k} Q_1^{m} Q_2^{n} Q_3^{p} Q_{123}^{q},\qquad
 \qquad i,j,m,n,p,q\in\mathbb{N},\qquad k\in\{0,1\},
\end{equation}
are linearly independent in $U_{\alpha}(\sl_2)^{\hat{\otimes} 3}$. First it is more
convenient (and equivalent) to replace the generators $Q_I$ by modified versions
analoguous to what was done in Section \ref{sec:szh}:
\begin{align}\label{}
 K_{I}:=\frac{Q_I-(q+q^{-1})}{(q-q^{-1})^{2}},\qquad I \in \{ 1,2,3,123,12,23\}.
\end{align}
The index $13$ does not belong to this set, and for this one, we set:
\begin{align}\label{}
 K_{13}:=[K_{12},K_{23}]_q.
\end{align}
Calculating explicitly the first terms in the expansions in $\alpha$,
we find that the new elements $K_{I}$ are
well-defined in  $U_{\alpha}(\sl_2)^{\hat{\otimes} 3}$, and moreover that their degree $0$
coefficients are expressed in terms of polarized traces, using the notations in
(\ref{gen-Racah}), as follows
\begin{gather}\label{}
 K_i|_{\alpha=0}=\frac{1}{2}k_i\quad (i=1,2,3),\qquad
 K_{123}|_{\alpha=0}=\frac{1}{2}k_4, \\
 K_{12}|_{\alpha=0}=\frac{1}{2}X,\qquad K_{23}|_{\alpha=0}=\frac{1}{2}Y,\qquad
 K_{13}|_{\alpha=0}=\frac{1}{4}Z.
\end{gather}
Now, to prove that the elements of the set (\ref{setQ}), with $Q_I$ replaced by $K_I$, are
linearly independent in $U_{\alpha}(\sl_2)^{\hat{\otimes} 3}$, it is enough to prove that
their ``classical limits'' (the degree 0 coefficients) are linearly independent in
$U(\sl_2)^{\otimes 3}$.
In view of the above calculations, it remains to show that the following set:
\begin{equation}\label{setT}
 k_1^ik_2^jk_3^kk_4^mX^nY^pZ^q, \qquad\qquad i,j,k,m,n,p\in\mathbb{N},\qquad q\in\{0,1\},
\end{equation}
is linearly independent in $U(\sl_2)^{\otimes 3}$.

\subsection{Racah algebra and diagonal centraliser in $U(\sl_2)^{\otimes 3}$}

To prove that the set (\ref{setT}) is linearly independent, we use the same line of
arguments as the one used in the study of the recoupling of two copies of $\sl(3)$. Thus
we only give here a sketch and refer for more details to
\cite{CrampePoulaindAndecyetal2020}.

It is known from classical invariant theory \cite{BereleStembridge1999, Drensky2006} that
the centralizer of the diagonal embedding of $U(\sl_2)$ in $U(\sl_2)^{\otimes 3}$ is
generated by the polarised traces $T^{(i,i)}$, $T^{(k,l)}$, $T^{(1,2,3)}$, with $i=1,2,3$
and $1\leq k<l\leq 3$, and moreover that the Hilbert--Poincar\'e series of the centralizer
is:
\begin{align}\label{eq:hilbert_centralizer}
 \frac{1-t^6}{(1-t^2)^6(1-t^3)}.
\end{align}
This series records the dimension for each degree of the centralizer, where the degree in
$U(\sl_2)^{\otimes3}$ is defined by $deg(e_{ij}^{(a)}) = 1$.
From this information, we extract at once that the set $k_1$, $k_2$, $k_3$, $k_4$, $X$,
$Y$, $Z$ generates the centralizer. Now, we have that these elements satisfy the classical
Racah relations:
\begin{equation}\label{relR1}
\begin{aligned}
 k_1,k_2,&k_3,k_4\quad \text{commute with all generators},\\
 [X,Y]&=Z,\\
 [X,Z]&=4\{X,Y\}+4X^2-4(k_1+k_2+k_3+k_4)X+4(k_1-k_2)(k_4-k_3),\\
 [Z,Y]&=4\{X,Y\}+4Y^2-4(k_1+k_2+k_3+k_4)Y+4(k_3-k_2)(k_4-k_1),
\end{aligned}
\end{equation}
together with
\begin{equation}\label{relR2}
 \Gamma=8(k_1-k_2+k_3-k_4)(k_1k_3-k_2k_4)-32(k_1k_3+k_2k_4),
\end{equation}
where the element $\Gamma$ is
\begin{align}
\begin{aligned}
 \Gamma:&=Z^2-8(XYX+YXY)+4(k_1+k_2+k_3+k_4-4)\{X,Y\}\\
 &~-8(k_1-k_2)(k_4-k_3)Y-8(k_3-k_2)(k_4-k_1)X.
\end{aligned}
\end{align}
The relations (\ref{relR1}) allow to rewrite any product in terms of ordered monomials in
the generators and (\ref{relR2}) allows to rewrite $Z^2$. So we deduce easily that the set
(\ref{setT}) is a spanning set for the centralizer. Finally, the comparison with the
Hilbert--Poincar\'e series in (\ref{eq:hilbert_centralizer}) shows that this set must be
linearly independent.

This concludes the proof of the injectivity of the map from $\saw(3)$ to $\mC_3\subset
U_{\alpha}(\sl_2)^{\hat{\otimes} 3}$.

\begin{rem}
Specializing the central elements $k_i$ to $\frac{m_i^2-1}{2}$, one finds that the
relations (\ref{relR1})--(\ref{relR2}) are expressed in terms of the polynomials:
\begin{align}\label{}
 \sum_{i=1}^{4}m_i^2,\qquad \left\{
 \begin{array}{c}
 \!\!(m_1^2-m_2^2)(m_4^2-m_3^2),\\[0.5em]
 \!\!(m_3^2-m_2^2)(m_4^2-m_1^2),
 \end{array}\right.\qquad
 (m_1^2m_3^2-m_2^2m_4^2)(m_1^2-m_2^2+m_3^2-m_4^2).
\end{align}
These polynomials are invariant polynomials under the action of the Weyl group $W(D_4)$ of
Section \ref{sec:szh}. This recovers explicitly the classical limit of the results in
Section \ref{sec:szh}.
\end{rem}

\printbibliography
\end{document}